\documentclass[12pt,twoside]{journal}
\usepackage{hyperref}
\usepackage{amsmath,amssymb,amsfonts} 
\usepackage{amsthm}
\usepackage{titlesec}
\usepackage[margin=1in]{geometry}
\usepackage{epsfig}
\usepackage{subfigure}
\numberwithin{equation}{section}
\usepackage{url}
\usepackage{enumerate}

\usepackage{color,colordvi}
\usepackage{psfrag}
\usepackage{color}
\usepackage{colordvi}
\usepackage{authblk}

\usepackage{fancyhdr}
\newcommand{\bsym}[1]{\boldsymbol{#1}}

\titleformat{\section}                                     
  {\sffamily\bfseries}
  {\thesection.}{1em}{}
\titleformat{\subsection}                                   
  {\sffamily\slshape\bfseries}
  {\thesubsection.}{1em}{}

\renewcommand\footnotemark{} 
\author[1,*]{L. L.~Yaw \thanks{*Correspondence to:  L. L. Yaw, Engineering Department, Walla Walla University, 100 SW 4th St, College Place, WA, 99324 USA.  E-mail: \texttt{\href{mailto:louie.yaw@wallawalla.edu}{louie.yaw@wallawalla.edu}}}}

\affil[1]{Engineering Department, Walla Walla University, 100 SW 4th St, College Place, WA 99324, USA}

\date{ }
\title{Introduction to the Virtual Element Method for 2D Elasticity} 

\begin{document}
\maketitle

\pagenumbering{arabic}

\section*{Summary} An introductory exposition of the virtual element method (VEM) is provided.  The intent is to make this method more accessible to those unfamiliar with VEM.  Familiarity with the finite element method for solving 2D linear elasticity problems is assumed.  Derivations relevant to successful implementation are covered.  Some theory is covered, but the focus here is on implementation and results.  Examples are given that illustrate the utility of the method.  Numerical results are provided to help researchers implement and verify their own results.\\

\noindent
KEY WORDS:  virtual element method, VEM, consistency, stability, polynomial base, polygon, vertices, elasticity, polymesher\\

\section{Introduction}

The virtual element method (VEM) originated around 2013 \cite{Beirao}.  It is yet another numerical method to solve partial differential equations.  VEM has many similarities with the finite element method (FEM).  One key difference is that VEM allows the problem domain to be discretized by a collection of arbitrary polygons.  The polygons need not all have the same number of sides.  Hence, one can have triangles, quadrilaterals, pentagons, and so on.  This is attractive as it makes meshing the problem domain easier using, for example, a Voronoi tesselation.  Convex and concave polygons are allowed. Linear, quadratic, and higher polynomial consistency is allowed within the method if implemented.  In this introductory exposition the focus is on solving 2D linear elasticity using linear order ($k=1$) polynomial interpolation.  Additionally, VEM has the ability to handle non-conforming discretizations (see Mengolini \cite{mengolini}).  For this document the goal is to provide implementation details and example results.  An attempt is made to present the information in a logical and meaningful order and thus provide the rationale for the method.  However, it is almost certainly not the order in which the method was discovered or rationalized originally.  For attributes not covered the interested reader is referred to the provided references.  Derivations and notation closely follows the paper by Mengolini et al. \cite{mengolini}, with some exceptions.

\section{The Continuous 2D Linear Elasticity Problem}
The goal is to solve 2D elasticity problems using VEM.
The weak form for the elasticity problem is:  Find $\mathbf{u} \in \boldsymbol{\mathcal{V}}$ such that
\begin{equation}\label{E1}
	a(\mathbf{u},\mathbf{v})=L(\mathbf{v}) \quad \forall \mathbf{v}\in \boldsymbol{\mathcal{V}},
\end{equation}
where the bilinear form is
\begin{equation}\label{E2}
	a(\mathbf{u},\mathbf{v})=\int _{\Omega} \boldsymbol{\sigma}(\mathbf{u}):\boldsymbol{\epsilon}(\mathbf{v})d\Omega,
\end{equation}
and the linear form is
\begin{equation}\label{E3}
	L(\mathbf{v})=\int_{\Omega} \mathbf{v}\cdot\mathbf{f}d\Omega+\int_{\partial\Omega_t} \mathbf{v}\cdot\bar{\mathbf{t}}d\partial\Omega.
\end{equation}\\

\noindent
\textbf{Remarks}
\begin{enumerate}[(i)]
	\item The vector-valued function space $\boldsymbol{\mathcal{V}}$ has components $v_1$ and $v_2$ that belong to the first-order Sobolev space $\mathcal{H}^1(\Omega)$ with zero values on displacement boundaries.
	\item The function $\mathbf{u} \in \boldsymbol{\mathcal{V}}$ is a trial solution, $\mathbf{v} \in \boldsymbol{\mathcal{V}}$ is a weight function.
\end{enumerate}

\section{Discretization of the problem domain}
\begin{figure}
\centering
\subfigure[]{\epsfig{file =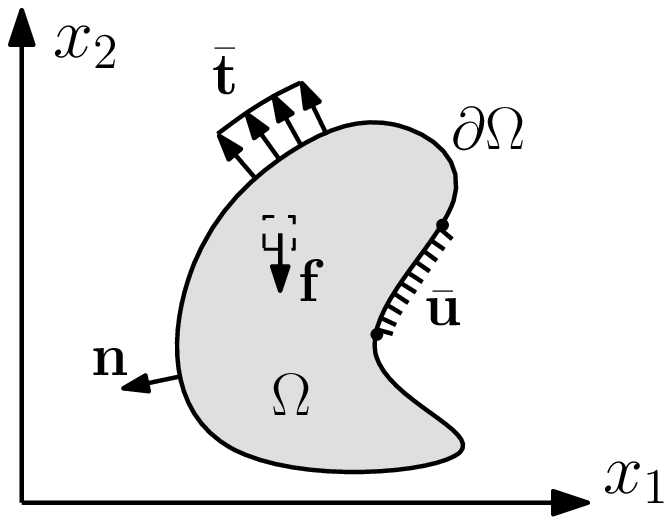,width=0.4\textwidth}}
\subfigure[]{\epsfig{file = 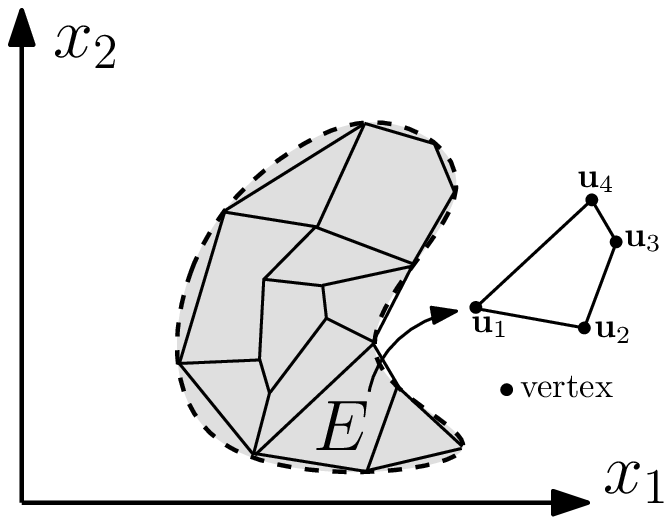,width=0.4\textwidth}}
\caption{2D Solid Domain: (a)~Elasticity problem with boundary conditions, (b)~Virtual element method domain discretization and example polygonal element with vector of nodal displacements labeling each vertex.}
\label{fig1}
\end{figure}

For 2D elasticity the domain is the geometric region, Figure~\ref{fig1}a, for which stresses, strains, and displacements are calculated.  In VEM, the domain is discretized with an arbitrary number of polygons, Figure~\ref{fig1}b.  Unlike FEM, polygon elements, convex or non-convex, with an arbitrary number of sides are used in VEM.  Due to the expectation that arbitrary polygon elements are used, it is necessary to imagine a space of interpolation functions that include polynomials but may also include non-polynomial functions.  This is necessary because the polygon elements must interconnect compatibly along their sides.  Hence, along the edges the interpolation functions are polynomials, but on the polygon interior the functions are possibly non-polynomial.  It turns out that it is not necessary to know the interpolation functions on the interior of the polygon elements, rather it is sufficient to know the polynomial functions along the polygon edges only.  The order of the polynomials along the edges are chosen at the beginning of the formulation.  As already indicated, this document chooses first order polynomials.
		
\section{VEM Functions}
The \emph{discrete} space of VEM functions,  $\boldsymbol{\mathcal{V}}^h$, over individual elements are a subset of the space of functions, $\boldsymbol{\mathcal{V}}$. The functions contained in $\boldsymbol{\mathcal{V}}$ can satisfy the weak form of the continuous 2D elasticity problem.  The superscript $h$ indicates a space of functions used as part of a discretization.  Mathematically, $\boldsymbol{\mathcal{V}}^h \subset \boldsymbol{\mathcal{V}}$.

A typical VEM function $\mathbf{v}^h \in \boldsymbol{\mathcal{V}}^h$ in 2D is a vector-valued \emph{displacement} function of spatial dimensions $(x_1,x_2)$.  For example, $\mathbf{v}^h=[v_1(x_1,x_2), v_2(x_1,x_2)]$ is a two component vector.  To represent such a function, basis (or shape) functions are needed for each element of the discretization.  A basis for the VEM space of functions within an element, along \emph{one} spatial dimension, is represented as $\{\varphi_i\}_{i=1,...,n_d}$, with $n_d$ degrees of freedom.  To organize this more clearly along \emph{both} spatial directions (2D case) a vector-valued form for the basis is written as $\bsym{\varphi}_1=[\varphi_1, 0]$, $\bsym{\varphi}_2=[0, \varphi_1]$,..., $\bsym{\varphi}_{2i-1}=[\varphi_i, 0]$, $\bsym{\varphi}_{2i}=[0, \varphi_i]$,..., $\bsym{\varphi}_{2n_d-1}=[\varphi_{n_d}, 0]$, $\bsym{\varphi}_{2n_d}=[0, \varphi_{n_d}]$.  Consequently, a VEM displacement function within an element written in terms of the basis functions is
\begin{equation}\label{v1}
	\mathbf{v}^h=\sum_{j=1}^{2n_d} \ dof_j(\mathbf{v}^h) \bsym{\varphi}_j 
\end{equation}
where here the operator $dof_j$ extracts the value of $\mathbf{v}^h$ at the $j$th degree of freedom.  As expected \eqref{v1} is a linear combination of basis functions.\\

\noindent
\textbf{Remarks}
\begin{enumerate}[(i)]
	\item VEM basis (shape) functions have the following characteristics: 
		\begin{itemize}
			\item continuous polynomial components of degree $k$ along element edges, 
			\item composed of polynomial functions and possibly non-polynomial functions on the interior of the element, this is why they are said to be unknown, and is why the word \emph{virtual} is used in the name of the method, 
			\item square integrable up to and including first derivatives, 
			\item Laplacian, $\Delta\mathbf{v}^h|_E$, is made of polynomials of degree $k-2$ in the interior of element~$E$, 			
			\item Kronecker delta property
		\end{itemize}
	\item VEM basis functions on the interior of the element can be found numerically by solving a PDE $\Delta\mathbf{v}^h=f$, with $f$ being a polynomial and prescribing values along the element's boundary.  This is not necessary, is costly, and is avoided to accomplish VEM. 
	\item $dof_i(\varphi_j)=\delta_{ij}$, this is enforced by how assumed characteristics of the basis functions are implemented in the calculations along with the operator $dof_i$.
	\item The $dof$ operator can refer to different components of the function, in the case of vector valued functions.
	\item The equation \eqref{v1} is similar to how interpolation between nodal values is accomplished in FEM.  Importantly, the basis functions are not actually known.  Although \eqref{v1} is a familiar form, it cannot be used.  Instead the VEM functions are projected onto a space of polynomial functions with a projection operator.  This is done with appropriate restrictions and adjustments in place to account for the fact that the 'correct' VEM functions aren't being used directly.
	\item In this document, element degrees of freedom are only considered at element vertices.  Hence, for 2D elasticity the vertices (or nodes) have 2 degrees of freedom (one in each coordinate direction $(x_1, x_2)$).  This is due to the choice of only using first order polynomials along the element boundaries.  More degrees of freedom per element and higher order polynomials are possible (see \cite{mengolini}).
\end{enumerate}

\section{Polynomial Functions}
With the concept of VEM functions realized, but knowing that they are not actually in hand, a clever strategy is to imagine the projection of the VEM functions onto a space of polynomial functions.  As it unfolds, in later sections, the strategy proves to be useful. First, it is useful, since a polynomial basis for the space of polynomial functions is easily created.  Second, because a very specific condition allows a projection operator to be found. Last, because conforming polynomials along the boundary are in line with VEM function assumed behavior and experience with FEM interpolation functions.

A space of scalar-valued polynomials of order equal to $k$ or less on an element $E$ is denoted as $\mathcal{P}_k(E)$.  This is extended to a 2D vector space of polynomials in two variables $\boldsymbol{\mathcal{P}}_k\equiv [\mathcal{P}_k]^2$. The polynomial space has a basis $\mathbf{P}_k=\{\mathbf{p}_{\alpha}\}_{\alpha=1,...,n_k}$.  An example case for polynomials of order $k=1$, clarifies the meaning.
\begin{equation}\label{pf1}
	\mathbf{P}_1=[\mathbf{p}_1,\ \mathbf{p}_2,\ \mathbf{p}_3,\ \mathbf{p}_4,\ \mathbf{p}_5,\ \mathbf{p}_6]
\end{equation}
or
\begin{equation}\label{pf2}
	\mathbf{P}_1=\left[\left ( \begin{array}{c}
	1\\
	0\\
	\end{array}\right ),
	\ \left ( \begin{array}{c}
	0\\
	1\\
	\end{array}\right ),
	\ \left ( \begin{array}{c}
	-\eta\\
	\xi\\
	\end{array}\right ),
	\ \left ( \begin{array}{c}
	\eta\\
	\xi\\
	\end{array}\right ),
	\ \left ( \begin{array}{c}
	\xi\\
	0\\
	\end{array}\right ),
	\ \left ( \begin{array}{c}
	0\\
	\eta\\
	\end{array}\right )
	\right].
\end{equation}
In the preceding equations, scaled monomials are used to construct the components of the polynomial basis of order $k=1$.  They are defined as
\begin{equation}
	\xi=\left(\frac{x_1-\bar{x}_1}{h_E}\right), \quad \eta=\left(\frac{x_2-\bar{x}_2}{h_E}\right),
\end{equation}
where $\bar{\mathbf{x}}=(\bar{x}_1,\bar{x}_2)$ is the centroid location of  element $E$ and $h_E$ is its diameter (i.e., diameter of smallest circle that encloses all vertices of the element).\\

\noindent
\textbf{Remarks}
\begin{enumerate}[(i)]
	\item In this document only polynomials of order $k=1$ are considered.  Higher order polynomials are possible (see \cite{mengolini}).
	\item With a choice of degrees of freedom the polynomials are unambiguously defined.  For example, two points make a line (a first order polynomial).  That is, nodal values are at vertices and first order polynomials interpolate between vertices of a given element edge.
	\item Just like $\mathbb{R}^2$ represents the space of 2D vectors, $[\mathcal{P}_k]^2$ represents a 2D vector polynomial with components in two variables.
	\item The number of terms (cardinality) in a polynomial base is calculated as $n_k=(k+1)(k+2)$.  For the case of $k=1$, $n_k=6$, which matches the number of terms in the polynomial base $\mathbf{P}_1$.  This arises by the requirement that all monomials of Pascal's triangle of order less than or equal to $k$ are included.
	\item Infinitesimal rigid body motions are represented in the first three monomials $\mathbf{p}_1, \mathbf{p}_2, \mathbf{p}_3$ of \eqref{pf2}.
	\item Based on $\alpha=1,2,3$ of the polynomial base the infinitesimal strain equals zero, $\bsym{\epsilon}(\mathbf{p}_{\alpha})=~\mathbf{0}$, since these terms are associated with rigid body motion.  To see this, recall that the strain displacement relations are often represented as $\bsym{\epsilon}=\mathbf{B}\mathbf{u}^E=\bsym{\partial}\mathbf{N}\mathbf{u}^E$.  Here, $\mathbf{v}^h=[v_1 \ v_2]^T$ is used as the displacement vector and $\mathbf{u}^E$ as the nodal(vertex) values of a typical polygon element. Then (with Voigt notation in mind) the differential operator, vector of shape functions, strain displacement matrix, and  vector of nodal values are as follows:
	
	\begin{equation}\label{pf3}	
			\bsym{\partial}=\left [ \begin{array}{cc}
		\partial_{x_1} & 0\\
		0 & \partial_{x_2}\\
		\partial_{x_2} & \partial_{x_1}\\
		\end{array}\right ]
	\end{equation}
	
		\begin{equation}\label{pf4}	
			\mathbf{N}=\left [ \begin{array}{cccc}
			\varphi_1 & \varphi_2 & ... & \varphi_{n_d}\\
			\end{array}\right ]	
		\end{equation}
	
		\begin{equation}\label{pf5}	
			\mathbf{B}=\bsym{\partial}\mathbf{N} = \left [ \begin{array}{cccc}
		\bsym{\partial}\varphi_1 & \bsym{\partial}\varphi_2 & ... & \bsym{\partial}\varphi_{n_d}\\
		\end{array}\right ]
	\end{equation}

	\begin{equation}\label{pf6}	
			\mathbf{u}^E=\left[ \begin{array}{ccccc}
		u^1_1 & u^1_2 & ... & u^{2n_d}_1 & u^{2n_d}_2
		\end{array}\right]^T	
	\end{equation}	
	
	Finally, with the above in hand, the engineering strains are written with the strain operator as
	
		\begin{equation}\label{pf7}	
			\bsym{\epsilon}=\left [ \begin{array}{c}
		\partial_{x_1}v_1\\
		\partial_{x_2}v_2\\
		\partial_{x_2}v_1 + \partial_{x_1}v_2\\
		\end{array}\right ]
		=\mathbf{B}\mathbf{u}^E
		=\left [ \bsym{\partial}\mathbf{N}\right ]\mathbf{u}^E = \left [ \begin{array}{cccc}
		\bsym{\partial}\varphi_1 & \bsym{\partial}\varphi_2 & ... & \bsym{\partial}\varphi_{n_d}\\
		\end{array}\right ]\mathbf{u}^E
	\end{equation}
	\item It is important to note that in the preceding equations, \eqref{pf3} to \eqref{pf7}, VEM basis functions are used conceptually.  Yet, these functions are not known.  In fact, it is necessary to insert the projection of VEM basis functions.  In later sections this is discussed further.  Toward this end, the projection operator is determined next.	
\end{enumerate}

\section{The Projector}
A descretization using polygons is the goal.  Functions that fit the necessary conditions on polygons are called VEM functions.  These functions are not known in a form that allows implementation in a numerical formulation.  However, it is possible to recover an approximation of the VEM functions by projecting them onto a polynomial basis.  Importantly, the VEM functions projected onto the polynomial basis create polynomials along the element edges and are able to exactly reproduce polynomials up to order $k$.  This is called k-\emph{consistency}.  Consistency and stability are both required for the success of a numerical discretization.  Stability is addressed later.  Nevertheless, to achieve consistency and for the projection to provide the best approximation of the VEM functions, the following orthogonality criteria using the projection operator, $\Pi$, in each polygon is enforced:

\begin{equation}\label{p1}
		a_E(\mathbf{u}^h-\Pi\mathbf{u}^h,\mathbf{p})=0, \quad \forall \mathbf{p} \in \bsym{\mathcal{P}}_k(E),
\end{equation}
where trial solution $\mathbf{u}^h \in \bsym{\mathcal{V}}^h$.
Recall the bilinear form \eqref{E2}. The terms in $a_E$ are inserted in \eqref{E2} and the integration takes place over an individual element $E$.  The result is a measure of strain energy.  In \eqref{p1} $\mathbf{u}^h-\Pi\mathbf{u}^h$ is the error (or difference) between the VEM function and the projection.  In the ensuing derivations the projector is solved for by using \eqref{p1} so that the error is orthogonal to each polynomial basis in the polynomial space $\bsym{\mathcal{P}}_k(E)$.  In essence, this implies that the energy error is not captured by the polynomial basis.  In other words, the polynomial basis is forced to not include any of the energy error caused by using the projection.  This is exactly how k-consistency is enforced.  The equation \eqref{p1} is the starting point for finding the projector matrix.  Furthermore, unless explicitly stated, the simplification $\Pi\equiv\Pi_{E,k}$ is implied.  The symbol $\Pi_{E,k}$ projects element functions from the VEM space onto the space of polynomials of order $k$, mathematically, $\Pi_{E,k}: \bsym{\mathcal{V}}^h(E)\rightarrow \bsym{\mathcal{P}}_k(E)$.

To solve for the projector begin by rearranging \eqref{p1}.
\begin{equation}\label{p2}
	\begin{split}
	&a_E(\mathbf{u}^h,\mathbf{p})-a_E(\Pi\mathbf{u}^h,\mathbf{p})=0\\
	\Rightarrow \quad &a_E(\mathbf{u}^h,\mathbf{p})=a_E(\Pi\mathbf{u}^h,\mathbf{p}), \quad \forall \mathbf{p} \in \bsym{\mathcal{P}}_k(E).
	\end{split}
\end{equation}

Then substituting terms into \eqref{p2}, and noting nodal displacements cancel from both sides and that the strain operator is linear, yields
\begin{equation}\label{p3}
	\int _{E} \bsym{\epsilon}(\bsym{\varphi}_i)^T\mathbf{C}\bsym{\epsilon}(\mathbf{p}_{\alpha}) dE
	=\int _{E} \bsym{\epsilon}(\Pi(\bsym{\varphi}_i))^T\mathbf{C}\bsym{\epsilon}(\mathbf{p}_{\alpha}) dE.
\end{equation}
Since the projection is onto the space of polynomials, it is reasonable to replace it with a linear combination of polynomial basis functions.  To this end, observe
\begin{equation}\label{p4}
	\Pi(\bsym{\varphi}_i)=\sum\limits^{n_k}_{\beta=1}s_{i,\beta}\mathbf{p}_{\beta} \quad \quad i=1,...,2n_d.
\end{equation}
Inserting \eqref{p4} into \eqref{p3} yields
\begin{equation}\label{p5}
	\int _{E} \bsym{\epsilon}(\bsym{\varphi}_i)^T\mathbf{C}\bsym{\epsilon}(\mathbf{p}_{\alpha}) dE
	=\sum\limits^{n_k}_{\beta=1}s_{i,\beta}\int _{E} \bsym{\epsilon}(\mathbf{p}_{\beta})^T\mathbf{C}\bsym{\epsilon}(\mathbf{p}_{\alpha}) dE.
\end{equation}
The above equation, for a particular value of VEM shape function $i$, gives $\alpha=1,...,n_k$ simultaneous linear equations with $n_k$ unknowns $s_{i,\beta}$.  This is written as
\begin{equation}\label{p6}
	b_{i,\alpha}=\sum\limits^{n_k}_{\beta=1}s_{i,\beta}\tilde{G}_{\alpha\beta}.
\end{equation}
In matrix form \eqref{p6} becomes
\begin{equation}\label{p7}
	\mathbf{b}_i=\tilde{\mathbf{G}}\mathbf{s}_i,
\end{equation}
where
\begin{equation}\label{p8}
	\mathbf{b}_i=\left [ \begin{array}{c}
		a_E(\mathbf{p}_1,\bsym{\varphi}_i)\\
		\vdots \\
		a_E(\mathbf{p}_{n_k},\bsym{\varphi}_i)\\
		\end{array}\right ],\quad
		\mathbf{s}_i=\left [ \begin{array}{c}
		s_{i,1}\\
		\vdots \\
		s_{i,n_k}\\
		\end{array}\right ]
\end{equation}
and
\begin{equation}\label{p9}
	\tilde{G}_{\alpha\beta}=\int _{E} \bsym{\epsilon}(\mathbf{p}_{\beta})^T\mathbf{C}\bsym{\epsilon}(\mathbf{p}_{\alpha}) dE \quad \Rightarrow \quad \tilde{\mathbf{G}}=[n_k \ \times \ n_k], \quad \tilde{\mathbf{G}}=\tilde{\mathbf{G}}^T.
\end{equation}
Then recognizing that $i=1,...,2n_d$ the above equations are repeated for all values of $i$ so that
\begin{equation}\label{p10}
	\tilde{\mathbf{B}}= \tilde{\mathbf{G}}\tilde{\bsym{\Pi}}^*
\end{equation}

\begin{equation}\label{p11}
	\tilde{\mathbf{B}}= \left[ \begin{array}{cccc}
		\mathbf{b}_1 & \mathbf{b}_2& ... & \mathbf{b}_{2n_d} \\
		\end{array}\right], \quad \tilde{\mathbf{B}}=[n_k \times 2n_d]
\end{equation}

\begin{equation}\label{p12}
	\tilde{\bsym{\Pi}}^*= \left[ \begin{array}{cccc}
		\mathbf{s}_1 & \mathbf{s}_2& ... & \mathbf{s}_{2n_d}\\
		\end{array}\right], \quad \tilde{\bsym{\Pi}}^*=[n_k \times 2n_d].
\end{equation}
Yet, \eqref{p10} needs modification.  Observe, for $\alpha=1,2,3$ equation \eqref{p5} results in $0=0$.  This is because the strain terms evaluate to zero for the rigid body modes of the polynomial base.  Hence, \eqref{p10} is an undetermined system for $\tilde{\bsym{\Pi}}^*$.  Three additional equations are obtained by requiring that
\begin{equation}\label{p13}
	\begin{split}
		&\frac{1}{n_v}\sum\limits^{2n_v}_{i=1}dof_i(\mathbf{v}^h)dof_i(\mathbf{p}_{\alpha})=\frac{1}{n_v}\sum\limits^{2n_v}_{i=1}dof_i(\Pi(\mathbf{v}^h))dof_i(\mathbf{p}_{\alpha}), \quad \text{\normalfont for} \ \alpha=1,2,3\\
	\Rightarrow	&\frac{1}{n_v}\sum\limits^{2n_v}_{i=1}v_i^h dof_i(\bsym{\varphi}_I)dof_i(\mathbf{p}_{\alpha})=\frac{1}{n_v}\sum\limits^{2n_v}_{i=1}v_i^h dof_i\left(\sum\limits^{n_k}_{\beta=1}s_{I,\beta}\mathbf{p}_\beta\right)dof_i(\mathbf{p}_{\alpha})\\
	\Rightarrow	&\frac{1}{n_v}\sum\limits^{2n_v}_{i=1}dof_i(\bsym{\varphi}_I)dof_i(\mathbf{p}_{\alpha})=\frac{1}{n_v}\sum\limits^{2n_v}_{i=1}dof_i\left(\sum\limits^{n_k}_{\beta=1}s_{I,\beta}\mathbf{p}_\beta\right)dof_i(\mathbf{p}_{\alpha})\\
	\Rightarrow	&\frac{1}{n_v}\sum\limits^{2n_v}_{i=1}dof_i(\bsym{\varphi}_I)dof_i(\mathbf{p}_{\alpha})=\frac{1}{n_v}\sum\limits^{2n_v}_{i=1}\sum\limits^{n_k}_{\beta=1}s_{I,\beta}dof_i(\mathbf{p}_\beta)dof_i(\mathbf{p}_{\alpha})\\
	\end{split}
\end{equation}
The last line of \eqref{p13} in matrix form becomes
\begin{equation}\label{p14}
	\breve{\mathbf{b}}_I=\breve{\mathbf{G}}\breve{\mathbf{s}}_I,
\end{equation}
where
\begin{equation}\label{p15}
	\breve{\mathbf{b}}_I=\left [ \begin{array}{c}
		\frac{1}{n_v}\sum\limits^{2n_v}_{i=1}dof_i(\bsym{\varphi}_I)dof_i(\mathbf{p}_1)\\
		\frac{1}{n_v}\sum\limits^{2n_v}_{i=1}dof_i(\bsym{\varphi}_I)dof_i(\mathbf{p}_2)\\
		\frac{1}{n_v}\sum\limits^{2n_v}_{i=1}dof_i(\bsym{\varphi}_I)dof_i(\mathbf{p}_3)\\
		\end{array}\right ],\quad
		\breve{\mathbf{s}}_I=\left [ \begin{array}{c}
		\breve{s}_{I,1}\\
		\vdots \\
		\breve{s}_{I,n_k}\\
		\end{array}\right ],
\end{equation}

\begin{equation}\label{p16}
	\breve{\mathbf{B}}= \left[ \begin{array}{cccc}
		\breve{\mathbf{b}}_1 & \breve{\mathbf{b}}_2& ... & \breve{\mathbf{b}}_{2n_d} \\
		\end{array}\right], \quad \breve{\mathbf{B}}=[3 \times 2n_d],
\end{equation}

\begin{equation}\label{p17}
	\breve{\mathbf{\Pi}}= \left[ \begin{array}{cccc}
		\breve{\mathbf{s}}_1 & \breve{\mathbf{s}}_2& ... & \breve{\mathbf{s}}_{2n_d} \\
		\end{array}\right], \quad \breve{\mathbf{\Pi}}=[n_k \times 2n_d],
\end{equation}
and
\begin{equation}\label{p18}
	\breve{G}_{\alpha\beta}=\frac{1}{n_v}\sum\limits^{2n_v}_{i=1} dof_i(\mathbf{p}_\alpha)dof_i(\mathbf{p}_{\beta}) \quad \Rightarrow \quad \breve{\mathbf{G}}=[3 \ \times \ n_k].
\end{equation}
Consequently, the three equations in matrix form are
\begin{equation}\label{p19}
	\breve{\mathbf{B}}=\breve{\mathbf{G}}\breve{\mathbf{\Pi}}.
\end{equation}
Finally, \eqref{p10} is modified so that \eqref{p19} occupies the first three rows.  The final modified form of the equations is denoted as
\begin{equation}\label{p20}
	\bar{\mathbf{B}}=\mathbf{G}\tilde{\mathbf{\Pi}}.
\end{equation}
Hence, the projector is
\begin{equation}\label{p21}
	\tilde{\mathbf{\Pi}}=\mathbf{G}^{-1}\bar{\mathbf{B}}.
\end{equation}\\

\noindent
\textbf{Remarks}
\begin{enumerate}[(i)]
	\item The matrix $\bar{\mathbf{B}}$, defined in \eqref{p20}, should not be confused with the $\bar{\mathbf{B}}$ matrix \cite{Hughes} used in the finite element method for incompressibility problems.
	\item The terms in \eqref{p15} for $\breve{\mathbf{b}}_I$ simplify further to 	
	\begin{equation}\label{p22}		
		\breve{\mathbf{b}}_I=\left [ \begin{array}{c}
		\frac{1}{n_v}\sum\limits^{2n_v}_{i=1}\delta_{iI}dof_i(\mathbf{p}_1)\\
		\frac{1}{n_v}\sum\limits^{2n_v}_{i=1}\delta_{iI}dof_i(\mathbf{p}_2)\\
		\frac{1}{n_v}\sum\limits^{2n_v}_{i=1}\delta_{iI}dof_i(\mathbf{p}_3)\\
		\end{array}\right ], \quad \quad \delta_{iI}=\left\{ \begin{array}{c} 
		0 \ \text{for} \ i\neq I\\
		1 \ \text{for} \ i=I.
		\end{array}
		\right.
	\end{equation}	
\end{enumerate}

\section{Element Stiffness}
Recall the objective is to discretize the domain with polygon elements.  VEM functions with all the requisite characteristics are necessary to interpolate over the domain of each individual element.  The strain energy for an element is expressed in the discrete bilinear form as
\begin{equation}\label{e1}
	a_E(\mathbf{u}^h,\mathbf{v}^h)=\int _E \boldsymbol{\sigma}(\mathbf{u}^h):\boldsymbol{\epsilon}(\mathbf{v}^h)dE
	=\int _{E} \bsym{\epsilon}(\mathbf{v}^h)^T\mathbf{C}\bsym{\epsilon}(\mathbf{u}^h) dE.
\end{equation}
Yet the VEM functions are not known, and it is preferable to write the bilinear form in terms of the projection \cite{nguyen}.  With this motivation the bilinear form is written as

\begin{equation}\label{e2}
	\begin{split}
		a_E(\mathbf{u}^h,\mathbf{v}^h)&=a_E(\Pi\mathbf{u}^h+(\mathbf{u}^h-\Pi\mathbf{u}^h),\Pi\mathbf{v}^h+(\mathbf{v}^h-\Pi\mathbf{v}^h))\\
		&=a_E(\Pi\mathbf{u}^h,\Pi\mathbf{v}^h)+a_E(\mathbf{u}^h-\Pi\mathbf{u}^h,\Pi\mathbf{v}^h)\\
		&\quad \quad +a_E(\Pi\mathbf{u}^h,\mathbf{v}^h-\Pi\mathbf{v}^h)+a_E(\mathbf{u}^h-\Pi\mathbf{u}^h,\mathbf{v}^h-\Pi\mathbf{v}^h)\\
		&=\underbrace{a_E(\Pi\mathbf{u}^h,\Pi\mathbf{v}^h)}_{\text{\normalfont 1st part}}+\underbrace{a_E(\mathbf{u}^h-\Pi\mathbf{u}^h,\mathbf{v}^h-\Pi\mathbf{v}^h)}_{\text{\normalfont 2nd part}}.
	\end{split}
\end{equation}
In the last step of \eqref{e2} the other terms vanish due to the enforcement of equation \eqref{p1}.  

Partial success is now achieved in the last line of \eqref{e2}.  The first part is expressed entirely in terms of the projected VEM functions and leads to the consistent part of the stiffness matrix.  The second part leads to stiffness stability, which `corrects' for what is lost of the VEM functions due to the projection.  Each of the parts are dealt with in turn.

\subsection{Stiffness providing consistency}
It is possible to obtain the consistent part of the stiffness exactly since it is projected onto the known polynomial base functions. The first part of \eqref{e2}, similar to \eqref{e1}, leads to
\begin{equation}\label{e3}
	a_E(\Pi\mathbf{u}^h,\Pi\mathbf{v}^h)=\int _{E} \bsym{\epsilon}(\Pi\mathbf{v}^h)^T\mathbf{C}\bsym{\epsilon}(\Pi\mathbf{u}^h)dE.
\end{equation}
Then using \eqref{v1} and focusing on specific dofs $i$ and $j$
\begin{equation}\label{e4}
	a_E(\Pi\mathbf{u}^h,\Pi\mathbf{v}^h)_{i,j}=dof_i(\mathbf{v}^h)\underbrace{\int _{E} \bsym{\epsilon}(\Pi(\bsym{\varphi}_i))^T\mathbf{C}\bsym{\epsilon}(\Pi(\bsym{\varphi}_j))dE}_{(k^c_E)_{ij}} dof_j(\mathbf{u}^h).
\end{equation}
Now, taking the $ij$ component of the element stiffness from \eqref{e4}, equation \eqref{p4}, and linearity of the strain operator, observe
\begin{equation}\label{e5}
	\begin{split}
		(\mathbf{k}^c_E)_{ij}&=\int _{E} \bsym{\epsilon}(\Pi(\bsym{\varphi}_i))^T\mathbf{C}\bsym{\epsilon}(\Pi(\bsym{\varphi}_j))dE \\
		&=\int _{E} \bsym{\epsilon}\left(\sum\limits^{n_k}_{\alpha=1}s_{i,\alpha}\mathbf{p}_{\alpha}\right)^T\mathbf{C}\bsym{\epsilon}\left(\sum\limits^{n_k}_{\beta=1}s_{j,\beta}\mathbf{p}_{\beta}\right)dE \\
		&=\sum\limits^{n_k}_{\alpha=1}\sum\limits^{n_k}_{\beta=1}s_{i,\alpha}s_{j,\beta}\int _{E} \bsym{\epsilon}(\mathbf{p}_{\alpha})^T\mathbf{C}\bsym{\epsilon}(\mathbf{p}_{\beta})dE \\
		&=\sum\limits^{n_k}_{\alpha=1}\sum\limits^{n_k}_{\beta=1}s_{i,\alpha}s_{j,\beta}a_{E} (\mathbf{p}_{\alpha},\mathbf{p}_{\beta}) \\
		&=\sum\limits^{n_k}_{\alpha=1}\sum\limits^{n_k}_{\beta=1}\tilde{\bsym{\Pi}}_{\alpha,i}\tilde{\bsym{\Pi}}_{\beta,j}\tilde{\mathbf{G}}_{\alpha\beta} \\
		&=\left[\tilde{\bsym{\Pi}}^T\tilde{\mathbf{G}}\tilde{\bsym{\Pi}}\right]_{ij},		
	\end{split}	
\end{equation}
where
\begin{equation}
	\tilde{\mathbf{G}}=a_E(\mathbf{p}_{\alpha},\mathbf{p}_{\beta})=\int _{E} \bsym{\epsilon}(\mathbf{p}_{\alpha})^T\mathbf{C}\bsym{\epsilon}(\mathbf{p}_{\beta})dE.
\end{equation}
Consequently, the consistent part of the stiffness matrix is represented as
\begin{equation}\label{e6}
		\mathbf{k}^c_E=\tilde{\bsym{\Pi}}^T\tilde{\mathbf{G}}\tilde{\bsym{\Pi}}.
\end{equation}
 
\subsection{Stiffness providing stability}
To deal with the stability stiffness several constructions need to be set in place, motivated by the approach of Sukumar and Tupek \cite{suku:elastodyn}, yet with some matrices ordered to match the approach of Mengolini et al. \cite{mengolini}, as used herein.  First, define a matrix of VEM basis functions as
\begin{equation}\label{stab1}
	\bsym{\varphi}=\left[ \begin{array}{cccccccccc}
		\varphi_1 & 0          & \varphi_2 & 0         & \cdots & \varphi_i & 0         & \cdots & \varphi_{n_d} & 0\\
		0         & \varphi_1  & 0         & \varphi_2 & \cdots & 0         & \varphi_i & \cdots & 0             & \varphi_{n_d}   
		\end{array}\right]
		=[\bsym{\varphi}_1 \ \bsym{\varphi}_2 \cdots \bsym{\varphi}_{2n_d}].
\end{equation}
Then, considering \eqref{p4}, which relates the projector to the polynomial basis for a single basis function, $\bsym{\varphi}_i$, the corresponding matrix expression is
\begin{equation}\label{stab2}
	\Pi(\bsym{\varphi})=\Pi\{\bsym{\varphi}_1 \ \bsym{\varphi}_2 \ ... \ \bsym{\varphi}_{2n_d}\}=\sum\limits_{\beta=1}^{n_k}\mathbf{p}_{\beta}\{s_{1,\beta} \ s_{2,\beta} \ \cdots \ s_{2n_d,\beta} \}=\mathbf{P}_1\tilde{\bsym{\Pi}}.
\end{equation}
Equation \eqref{stab2} is an expression in matrix form that represents the projection of the VEM basis functions onto the polynomial basis.

Next, a $\mathbf{D}$ matrix is defined as
\begin{equation}\label{stab3}
	\mathbf{D}_{i\alpha}=dof_i(\mathbf{p}_{\alpha}).
\end{equation}
Observe that constant and linear reproducing conditions of the VEM basis provide the following relations:
\begin{equation}\label{stab4}
	\sum\limits_{i=1}^{n_d}\varphi_i(\mathbf{x})=1, \quad \sum\limits_{i=1}^{n_d}\varphi_i(\mathbf{x}) \xi_i=\xi, \quad \sum\limits_{i=1}^{n_d}\varphi_i(\mathbf{x}) \eta_i=\eta.
\end{equation}
The preceding reproducing conditions are used to express the polynomial base in terms of $\bsym{\varphi}$ and $\mathbf{D}$ as follows: 
\begin{equation}\label{stab5}
	\mathbf{P}_1=\bsym{\varphi}\mathbf{D}.
\end{equation}
To see how \eqref{stab5} comes about, it is instructive to write out the matrices $\bsym{\varphi}$ and $\mathbf{D}$ with internal components.  Observing how the components multiply together and sum, reveals the reproducing conditions. 
Finally, using \eqref{stab5} in \eqref{stab2} yields
\begin{equation}\label{stab6}
	\Pi(\bsym{\varphi})=\mathbf{P}_1\tilde{\bsym{\Pi}}=\bsym{\varphi}\mathbf{D}\tilde{\bsym{\Pi}}=\bsym{\varphi}\bsym{\Pi}.
\end{equation}
Note that equation \eqref{stab6} provides the matrix representation of the projection in two ways.  The first way is the projection onto the polynomial basis as found in \eqref{stab2}.  The second way is the projection on the $\bsym{\varphi}$ basis set, from which the projection matrix, $\bsym{\Pi}$ is defined as
\begin{equation}\label{stab7}
	\bsym{\Pi}=\mathbf{D}\tilde{\bsym{\Pi}}.
\end{equation}
This last form of the projection proves useful to determine the stability stiffness.

From the second part of \eqref{e2} and using \eqref{v1} it follows that
\begin{equation}\label{e7}
\begin{split}
	a_E&(\mathbf{u}^h-\Pi\mathbf{u}^h,\mathbf{v}^h-\Pi\mathbf{v}^h)_{ij}\\
	&=a_E(\bsym{\varphi}_j dof_j(\mathbf{u}^h)-\Pi(\bsym{\varphi}_j) dof_j(\mathbf{u}^h),\bsym{\varphi}_i dof_i(\mathbf{v}^h)-\Pi(\bsym{\varphi}_i) dof_i(\mathbf{v}^h))\\
	&=dof_i(\mathbf{v}^h)\underbrace{a_E(\bsym{\varphi}_j-\Pi(\bsym{\varphi}_j), \bsym{\varphi}_i-\Pi(\bsym{\varphi}_i))}_{(\mathbf{k}^s_E)_{ij}} dof_j(\mathbf{u}^h).
	\end{split}
\end{equation}
Hence, the stability part of the stiffness for dofs $i$ and $j$ is
\begin{equation}\label{e8}
	\begin{split}
	(\mathbf{k}^s_E)_{ij}&=a_E(\bsym{\varphi}_j-\Pi(\bsym{\varphi}_j), \bsym{\varphi}_i-\Pi(\bsym{\varphi}_i))\\
	&=a_E((1-\Pi)\bsym{\varphi}_j, (1-\Pi)\bsym{\varphi}_i).
	\end{split}
\end{equation}
In light of \eqref{stab6} and equation \eqref{e8} the complete stability stiffness is
\begin{equation}\label{e9}	
	\begin{split}
	\mathbf{k}^s_E&=a_E((1-\Pi)\bsym{\varphi}, (1-\Pi)\bsym{\varphi})\\
	&=a_E(\bsym{\varphi}(\mathbf{I}-\bsym{\Pi}), \bsym{\varphi}(\mathbf{I}-\bsym{\Pi})).	
	\end{split}
\end{equation}
Observing the similarity to \eqref{e5} the terms $(\mathbf{I}-\bsym{\Pi})$ are moved outside the bilinear form, so that \eqref{e9} becomes
\begin{equation}\label{e9a}
	\mathbf{k}^s_E=(\mathbf{I}-\bsym{\Pi})^Ta_E(\bsym{\varphi},\bsym{\varphi})(\mathbf{I}-\bsym{\Pi}).
\end{equation}
It is not possible to evaluate the term $a_E(\bsym{\varphi},\bsym{\varphi})$ because it contains VEM shape functions, which are not known.  Hence, the effect of this term is approximated \cite{Beirao} \cite{mengolini} \cite{Artioli} using the  scaling, $\tau^h \text{\normalfont tr}(\mathbf{k}_E^c)$, where $\tau^h$ is a user-defined parameter which is taken as 1/2 for linear elasticity.  The stability stiffness then is written as
\begin{equation}\label{e9b}
	\mathbf{k}^s_E=\tau^h \text{\normalfont tr}(\mathbf{k}_E^c)(\mathbf{I}-\bsym{\Pi})^T(\mathbf{I}-\bsym{\Pi}).
\end{equation}
In the above expression, \eqref{e9b}, $\mathbf{I}$ is the $2n_d$ by $2n_d$ identity matrix.

An alternative approach advocated by \cite{suku:elastodyn} is to approximate $a_E(\bsym{\varphi},\bsym{\varphi})$ with a $2n_d$ by $2n_d$ diagonal matrix, $\mathbf{S}^d_E$, scaled appropriately.  The terms along the diagonal are taken as: $(\mathbf{S}^d_E)_{ii}=\mathsf{max}(\alpha_0 \ \text{\normalfont tr}(\mathbf{C})/m, (\mathbf{k}^c_E)_{ii})$, where $m=3$ in 2D, $\mathbf{C}$ is the 2D modular matrix for plane strain or plain stress, $\text{\normalfont tr}$ denotes the trace operator, and $\alpha_0=1$ since the formulation uses scaled monomials associated with the elements whose diameters are on the order of 1.  With this in hand the stability stiffness is represented as
\begin{equation}\label{e11}
	\mathbf{k}^s_E=(\mathbf{I}-\bsym{\Pi})^T\mathbf{S}^d_E(\mathbf{I}-\bsym{\Pi}).
\end{equation}

\subsection{Final form of element stiffness}
With the stiffnesses in hand, the total stiffness for element $E$ is
\begin{equation}\label{e10}
	\mathbf{k}_E=\mathbf{k}^c_E+\mathbf{k}^s_E
\end{equation}

\subsection{Numerical Evaluation of Various Terms}
In this subsection derivations and details are provided to assist in the numerical implementation of VEM.  Useful simplified expressions are given.  In particular, this subsection focuses on specific matrices necessary to construct the VEM element stiffness matrices.  Other specific implementation details of a VEM computer program are discussed by Mengolini et al.~\cite{mengolini}.
\subsubsection{The $\tilde{\mathbf{B}}$ matrix}
Some discussion is necessary to illustrate how certain terms are calculated.  First, consider calculation of a typical term in the $\tilde{\mathbf{B}}$ matrix.  A typical term is (see \eqref{p5})

\begin{equation}\label{ne1}
	\tilde{B}_{\alpha i}=a_E(\boldsymbol{\varphi}_i,\mathbf{p}_{\alpha})=\underbrace{\int _{E} \bsym{\epsilon}(\bsym{\varphi}_i)^T\mathbf{C}\bsym{\epsilon}(\mathbf{p}_{\alpha}) dE}_{\text {\normalfont Voigt notation}}=\underbrace{\int _E \boldsymbol{\epsilon}(\bsym{\varphi}_i):\boldsymbol{\sigma}(\mathbf{p}_{\alpha})dE}_{\text{\normalfont tensor notation}}
\end{equation}
The symmetric part of $\nabla \bsym{\varphi}_i$ is $\bsym{\epsilon}(\bsym{\varphi}_i)$.  Then, since $\bsym{\sigma}$ is symmetric, it follows that
\begin{equation}\label{ne2}
	\nabla \bsym{\varphi}_i:\bsym{\sigma}=\bsym{\epsilon}(\bsym{\varphi}_i):\bsym{\sigma}.
\end{equation}
Therefore, \eqref{ne1} becomes (continuing with tensor notation)
\begin{equation}\label{ne3}
	\tilde{B}_{\alpha i}=\int _E \nabla \bsym{\varphi}_i:\boldsymbol{\sigma}(\mathbf{p}_{\alpha})dE.
\end{equation}
Next, observe that
\begin{equation}\label{ne4}
	\begin{split}
	\nabla (\bsym{\varphi}_i\cdot \boldsymbol{\sigma})&=\nabla \bsym{\varphi}_i:\boldsymbol{\sigma}+\bsym{\varphi}_i\cdot \nabla \boldsymbol{\sigma}\\
	\Rightarrow \nabla \bsym{\varphi}_i:\boldsymbol{\sigma} &=-\bsym{\varphi}_i\cdot \nabla \boldsymbol{\sigma}+\nabla (\bsym{\varphi}_i\cdot \boldsymbol{\sigma}) \\
	\Rightarrow \int _E \nabla \bsym{\varphi}_i:\boldsymbol{\sigma}dE &=-\int _E \bsym{\varphi}_i\cdot \nabla \boldsymbol{\sigma}dE+\int _E \nabla (\bsym{\varphi}_i\cdot \boldsymbol{\sigma})dE
	\end{split}
\end{equation}
Using the divergence theorem on the last line of \eqref{ne4} and substituting the result into \eqref{ne3}
\begin{equation}\label{ne5}
	\tilde{B}_{\alpha i}= -\int _E \bsym{\varphi}_i\cdot \nabla \boldsymbol{\sigma}(\mathbf{p}_{\alpha})dE+\int _{\partial E} \bsym{\varphi}_i\cdot \boldsymbol{\sigma}(\mathbf{p}_{\alpha})\mathbf{n}_e de, 
\end{equation}
where $\mathbf{n}_e$ is the outward unit normal to the element edge and $e$ denotes an element edge.  It is equation \eqref{ne5} that is numerically integrated to determine the entries in the $\tilde{\mathbf{B}}$ matrix.

Realize now that for $k=1$ only the boundary integral of \eqref{ne5} is nonzero.  As a result,
\begin{equation}\label{ne6}
	\tilde{B}_{\alpha i}=\int _{\partial E} \bsym{\varphi}_i\cdot \boldsymbol{\sigma}(\mathbf{p}_{\alpha})\mathbf{n}_e de. 
\end{equation}

	\begin{figure}
  \centering
  \epsfig{file=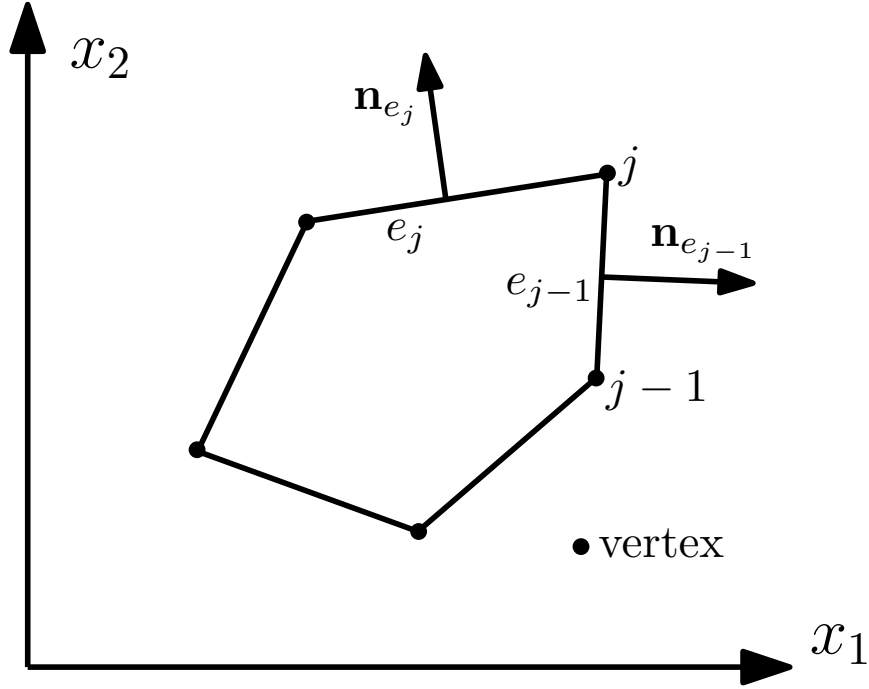,width=.7\textwidth}
  \caption{Single five sided element edges, normals, and nodes labeled.}\label{fige}
  \end{figure}  

Numerically, \eqref{ne6} is calculated by integrating around the boundary of the element (polygon) edges.  This is accomplished by using the vertex (node) points as the integration points and using the outward unit normal along each edge (see Figure \ref{fige}).  In essence, a trapezoidal rule is used to integrate along each polygon edge.  It is convenient to express the integration around the boundary as a sum over vertices
\begin{equation}\label{ne7}
	\tilde{B}_{\alpha i}=\sum\limits_{j=1}^{n_v} \bsym{\varphi}_i\cdot \boldsymbol{\sigma}(\mathbf{p}_{\alpha})\left(\frac{|e_{j-1}|}{2}\mathbf{n}_{e_{j-1}}+\frac{|e_{j}|}{2}\mathbf{n}_{e_{j}}\right). 
\end{equation}
where the stress terms are found by matrix multiplication
\begin{equation}\label{ne8}
	\begin{split}
		\boldsymbol{\sigma}(\mathbf{p}_{\alpha})&=\mathbf{C}\boldsymbol{\epsilon}(\mathbf{p}_{\alpha})\\
		\left[ \begin{array}{c}
		\sigma_x(\mathbf{p}_{\alpha})\\
		\sigma_y(\mathbf{p}_{\alpha})\\
		\sigma_{xy}(\mathbf{p}_{\alpha})
		\end{array}\right]&=\mathbf{C}
		\left[ \begin{array}{c}
		\epsilon_x(\mathbf{p}_{\alpha})\\
		\epsilon_y(\mathbf{p}_{\alpha})\\
		\gamma_{xy}(\mathbf{p}_{\alpha})
		\end{array}\right].
	\end{split}
\end{equation}
In the \eqref{ne8}, the appropriate $\mathbf{C}$ matrix for plane stress or plain strain is used (see equations \eqref{s2} and \eqref{s3}).  The result of \eqref{ne8} is then used to form the stress matrix
\begin{equation}\label{ne9}
	\boldsymbol{\sigma}(\mathbf{p}_{\alpha})=\left[ \begin{array}{cc}
		\sigma_x(\mathbf{p}_{\alpha}) & \sigma_{xy}(\mathbf{p}_{\alpha})\\
		\sigma_{xy}(\mathbf{p}_{\alpha}) &\sigma_y(\mathbf{p}_{\alpha})
		\end{array}\right].
\end{equation}
Then \eqref{ne9} is used in \eqref{ne7}
\begin{equation}\label{ne10}
	\tilde{B}_{\alpha i}=\sum\limits_{j=1}^{n_v} \bsym{\varphi}_i\cdot \left[ \begin{array}{cc}
		\sigma_x(\mathbf{p}_{\alpha}) & \sigma_{xy}(\mathbf{p}_{\alpha})\\
		\sigma_{xy}(\mathbf{p}_{\alpha}) &\sigma_y(\mathbf{p}_{\alpha})
		\end{array}\right]
		\left(\frac{|e_{j-1}|}{2}
		\left[ \begin{array}{c}
		n_{e1}\\
		n_{e2}
		\end{array}\right]_{j-1}+\frac{|e_{j}|}{2}
		\left[ \begin{array}{c}
		n_{e1}\\
		n_{e2}
		\end{array}\right]_j\right). 
\end{equation}

Finally, observe that $\bsym{\varphi}_i$ is nonzero only at node dofs $i=2j-1$ or $i=2j$, which correspond to two columns of the $\tilde{\mathbf{B}}$ matrix, so that 
\begin{equation}\label{ne11}
	\begin{split}
		\tilde{B}_{\alpha (2j-1)}&=\left[ \begin{array}{c}
		1\\
		0
		\end{array}\right]\cdot \left[ \begin{array}{cc}
		\sigma_x(\mathbf{p}_{\alpha}) & \sigma_{xy}(\mathbf{p}_{\alpha})\\
		\sigma_{xy}(\mathbf{p}_{\alpha}) &\sigma_y(\mathbf{p}_{\alpha})
		\end{array}\right]
		\left(\frac{|e_{j-1}|}{2}
		\left[ \begin{array}{c}
		n_{e1}\\
		n_{e2}
		\end{array}\right]_{j-1}+\frac{|e_{j}|}{2}
		\left[ \begin{array}{c}
		n_{e1}\\
		n_{e2}
		\end{array}\right]_j\right)\\ 
		\text{\normalfont and}&\\
		\tilde{B}_{\alpha (2j)}&=\left[ \begin{array}{c}
		0\\
		1
		\end{array}\right]\cdot \left[ \begin{array}{cc}
		\sigma_x(\mathbf{p}_{\alpha}) & \sigma_{xy}(\mathbf{p}_{\alpha})\\
		\sigma_{xy}(\mathbf{p}_{\alpha}) &\sigma_y(\mathbf{p}_{\alpha})
		\end{array}\right]
		\left(\frac{|e_{j-1}|}{2}
		\left[ \begin{array}{c}
		n_{e1}\\
		n_{e2}
		\end{array}\right]_{j-1}+\frac{|e_{j}|}{2}
		\left[ \begin{array}{c}
		n_{e1}\\
		n_{e2}
		\end{array}\right]_j\right).
	\end{split}
\end{equation}\\

\noindent
\textbf{Remarks}
\begin{enumerate}[(i)]
	\item Vertices $j$ range over $1$ to $n_v$ and \eqref{ne11} generates two columns of $\tilde{\mathbf{B}}$ at a time for the rows $\alpha=1$ to $n_k$.
	\item In other references the values at vertices are written in a slightly different form considering the normal to a line drawn between vertices $j-1$ and $j+1$.  However, for transparency the formula given above is provided.
	\item For vertex $j=1$ the edge length $|e_{j-1}|$ is taken as the the length between vertex $j=n_v$ and $j=1$, where $n_v$ is the number of vertices in the element (polygon).  Furthermore, the normal $\mathbf{n}_{e_{j-1}}$ is taken as the outward normal to the edge between $j=n_v$ and $j=1$.
	\item To be clear, the matrix multiplication of the right hand side of \eqref{ne11} results in a vector that is then dotted with either $[1 \ 0]^T$ or $[0 \ 1]^T$, as indicated.
	\item The $\tilde{\mathbf{B}}$ matrix is $n_k$ by $2n_v$ in size, for $k=1$.
	\item Equality \eqref{ne2} is true because, if $\mathbf{A}$ is an arbitrary tensor and $\mathbf{S}$ is a symmetric tensor, it can be shown that $\mathbf{A}:\mathbf{S}=\mathbf{A}^{sym}:\mathbf{S}$, where $\mathbf{A}^{sym}$ is the symmetric part of $\mathbf{A}$.
	\item In the above formulas, the more familiar subscripts $x,y$ for stresses and strains are used, which here refer to coordinate directions $x_1,x_2$, respectively.
\end{enumerate}

\subsubsection{The $\mathbf{D}$ matrix}
The $\mathbf{D}$ matrix is constructed by evaluating polynomial functions at the various degrees of freedom of polygon $E$.  The matrix entries are found by a straightforward evaluation of matrix terms.  The result is
\begin{equation}
	\mathbf{D}=\left[ \begin{array}{cccc}
		dof_1(\mathbf{p}_{1}) & dof_1(\mathbf{p}_{2})& \cdots & dof_1(\mathbf{p}_{n_k})\\
		dof_2(\mathbf{p}_{1}) & dof_2(\mathbf{p}_{2})& \cdots & dof_2(\mathbf{p}_{n_k})\\
		\vdots & \vdots & \ddots & \vdots\\
		dof_{2n_d}\mathbf{p}_{1}) & dof_{2n_d}(\mathbf{p}_{2})& \cdots & dof_{2n_d}(\mathbf{p}_{n_k})
		\end{array}\right].
\end{equation}

\subsubsection{The $\tilde{\mathbf{G}}$ matrix}
A typical term in the $\tilde{\mathbf{G}}$ matrix is expressed as
\begin{equation}\label{ne12}
	\tilde{G}_{\alpha \beta}=a_E(\mathbf{p}_{\alpha},\mathbf{p}_{\beta})=\underbrace{\int _{E} \bsym{\epsilon}(\mathbf{p}_{\beta})^T\mathbf{C}\bsym{\epsilon}(\mathbf{p}_{\alpha}) dE}_{\text {\normalfont Voigt notation}}=\underbrace{\int _E \boldsymbol{\epsilon}(\mathbf{p}_{\beta}):\boldsymbol{\sigma}(\mathbf{p}_{\alpha})dE}_{\text{\normalfont tensor notation}}.
\end{equation}
The previous work to find $\tilde{\mathbf{B}}$ is modified to find $\tilde{\mathbf{G}}$.  From \eqref{ne10} recognize that $\mathbf{p}_{\beta}$ in place of $\bsym{\varphi}_i$ yields
\begin{equation}\label{ne13}
	\tilde{G}_{\alpha \beta}=\sum\limits_{j=1}^{n_v} \mathbf{p}_{\beta}\cdot \left[ \begin{array}{cc}
		\sigma_x(\mathbf{p}_{\alpha}) & \sigma_{xy}(\mathbf{p}_{\alpha})\\
		\sigma_{xy}(\mathbf{p}_{\alpha}) &\sigma_y(\mathbf{p}_{\alpha})
		\end{array}\right]
		\left(\frac{|e_{j-1}|}{2}
		\left[ \begin{array}{c}
		n_{e1}\\
		n_{e2}
		\end{array}\right]_{j-1}+\frac{|e_{j}|}{2}
		\left[ \begin{array}{c}
		n_{e1}\\
		n_{e2}
		\end{array}\right]_j\right),
\end{equation}
where the quantities $\mathbf{p}_{\alpha}$ and $\mathbf{p}_{\beta}$ are evaluated at vertex $j$ for the $j$th term in the summation.  The matrix $\tilde{\mathbf{G}}$ is $n_k$ by $n_k$ in size.\\

\noindent
\textbf{Remarks}
\begin{enumerate}[(i)]
	\item It is possible to calculate the $\tilde{\mathbf{G}}$ matrix as $\tilde{\mathbf{G}}=\tilde{\mathbf{B}}\mathbf{D}$.  Hence, the above formulation for $\tilde{\mathbf{G}}$ provides an additional numerical check for verification.
	\item Observe that the first three rows of $\tilde{\mathbf{G}}$ and $\tilde{\mathbf{B}}$ need not be calculated because they contain all zeros.
	\item It is faster to calculate $\tilde{\mathbf{G}}$ by using $\tilde{\mathbf{G}}=\tilde{\mathbf{B}}\mathbf{D}$, once the algorithm is verified as working.  However, a better approach is to find $\mathbf{G}=\bar{\mathbf{B}}\mathbf{D}$ and then get $\tilde{\mathbf{G}}$ by zeroing the first three rows of $\mathbf{G}$.  The proof that $\mathbf{G}=\bar{\mathbf{B}}\mathbf{D}$ is shown in \cite{Beirao2} for one unknown per dof.  Here, a proof is given using the notation set forth so far and for 2D elasticity wherein two displacement unknowns per dof are present.\\
	
	Prove $\bar{\mathbf{B}}\mathbf{D}=\mathbf{G}$.
	\begin{proof}
	For $\alpha=1,2,3$, and making use of \eqref{p22}	
	\begin{equation}
		\sum_{i=1}^{2n_d}\bar{\mathbf{B}}_{\alpha i}\mathbf{D}_{i\beta}=\sum_{i=1}^{2n_d}\frac{1}{n_v}\delta_{ij}dof_j(\mathbf{p}_{\alpha})dof_i(\mathbf{p}_{\beta})=\sum_{i=1}^{2n_d}\frac{1}{n_v}dof_i(\mathbf{p}_{\alpha})dof_i(\mathbf{p}_{\beta})=\mathbf{G}_{\alpha \beta}.
	\end{equation}
	For $\alpha>3$	
	\begin{equation}
			\sum_{i=1}^{2n_d}\bar{\mathbf{B}}_{\alpha i}\mathbf{D}_{i\beta}=\sum_{i=1}^{2n_d}a_E(\mathbf{p}_{\alpha},\bsym{\varphi}_i)dof_i(\mathbf{p}_{\beta})=a_E(\mathbf{p}_{\alpha},\sum_{i=1}^{2n_d}dof_i(\mathbf{p}_{\beta})\bsym{\varphi}_i)=a_E(\mathbf{p}_{\alpha},\mathbf{p}_{\beta})=\mathbf{G}_{\alpha \beta}.
	\end{equation}
	Consequently, $\bar{\mathbf{B}}\mathbf{D}=\mathbf{G}$.
	\end{proof}
	
	\item To clarify, $\tilde{\mathbf{G}}$ is needed to obtain the stiffness matrix (see \eqref{e6}).  The terms $\mathbf{G}$ and $\bar{\mathbf{B}}$ are needed to calculate the projector in \eqref{p21}.
\end{enumerate}

\section{Application of External Forces}
External forces are caused by external tractions and body forces as indicated in equation \eqref{E3}.  Point loads are also possible.  All three forces are expressed for an individual element in the linear form 
\begin{equation}\label{ext1}
	L_E(\mathbf{v}^h)=\int_E \mathbf{v}^h\cdot\mathbf{f}dE+\int_{\partial E\cap\Omega_t} \mathbf{v}^h\cdot\bar{\mathbf{t}}d\partial E+\sum\limits_{i=1}\mathbf{v}^h(\mathbf{x}_i)\cdot\mathbf{F}_i. 
\end{equation}
The interested reader is directed to the discussion by Mengolini et al. \cite{mengolini}.  Herein, only external point loads are used in the examples shown in later sections.  From point loads a global external force vector is assembled, which is used to solve for the nodal displacements.  In a nonlinear analysis, the global external forces are used in a Newton-Raphson scheme to enforce equilibrium.

\section{Solving for Unknown Displacements}
Once element stiffness matrices are found they are assembled into a global stiffness matrix similar to FEM.  As a result the global stiffness is
\begin{equation}\label{sol1}
	\mathbf{K}=\overset{n_{elem}}{\underset{i=1}{\mathbf{\mathsf{A}}}}\mathbf{k}_E^{i}.
\end{equation}
Then with the global external force vector denoted as $\mathbf{F}$ the standard set of linear algebraic equations are
\begin{equation}\label{sol2}
	\mathbf{K}\mathbf{u}=\mathbf{F}.
\end{equation}
The global nodal displacements are then found in the typical manner
\begin{equation}\label{sol3}
	\mathbf{u}=\mathbf{K}^{-1}\mathbf{F}.
\end{equation}
	
\section{Element Strains}
Strains are found by starting with \eqref{v1}.  Then
\begin{equation}\label{ep1}
	\begin{split}
		\mathbf{v}^h\approx \Pi(\mathbf{v}^h)&=\sum_{j=1}^{2n_d} \ dof_j(\mathbf{v}^h) \Pi(\bsym{\varphi}_j)\\
		&=\left[ \begin{array}{cccc}
		\Pi(\bsym{\varphi}_1) & \Pi(\bsym{\varphi}_2)& ... & \Pi(\bsym{\varphi}_{2n_d}) \\
		\end{array}\right]\mathbf{u}^E,
	\end{split}
\end{equation}
where the vector of local values at dofs for the given element is
\begin{equation}\label{ep2}
	\mathbf{u}^E=\left[ \begin{array}{cccc}
		dof_1(\mathbf{v}^h) & dof_2(\mathbf{v}^h)& ... & dof_{2n_d}(\mathbf{v}^h)
		\end{array}\right]^T
		=\left[ \begin{array}{ccccc}
		u^1_1 & u^1_2 & ... & u^{2n_d}_1 & u^{2n_d}_2
		\end{array}\right]^T.
\end{equation}
Next the projection is expressed as
\begin{equation}\label{ep3}
	\Pi(\mathbf{v}^h)=\Pi(\bar{\mathbf{N}})\mathbf{u}^E
\end{equation}
where $\bar{\mathbf{N}}\equiv\bsym{\varphi}$ is the row vector of VEM basis functions
\begin{equation}\label{ep4}
	\bar{\mathbf{N}}=\left[ \begin{array}{cccc}
		\bsym{\varphi}_1 & \bsym{\varphi}_2 & ... & \bsym{\varphi}_{2n_d} \\
		\end{array}\right].
\end{equation}
Also, observe that \eqref{p4} leads to
\begin{equation}\label{ep5}
	\Pi(\bar{\mathbf{N}})=\left[ \begin{array}{cccc}
		\mathbf{p}_1 & \mathbf{p}_2 & ... & \mathbf{p}_{n_k} \\
		\end{array}\right]\tilde{\bsym{\Pi}}.
\end{equation}
Consequently,
\begin{equation}\label{ep6}
	\Pi(\mathbf{v}^h)=\left[ \begin{array}{cccc}
		\mathbf{p}_1 & \mathbf{p}_2 & ... & \mathbf{p}_{n_k} 
		\end{array}\right]\tilde{\bsym{\Pi}}\mathbf{u}^E.
\end{equation}
Last, using the strain operator \eqref{pf7}
\begin{equation}\label{ep7}
	\begin{split}
		\bsym{\epsilon}(\mathbf{v}^h)\approx\bsym{\epsilon}(\Pi(\mathbf{v}^h))&=\bsym{\epsilon}\left(\left[ \begin{array}{cccc}
		\mathbf{p}_1 & \mathbf{p}_2 & ... & \mathbf{p}_{n_k} 
		\end{array}\right]\tilde{\bsym{\Pi}}\mathbf{u}^E\right)\\
		&=\bsym{\epsilon}\left[ \begin{array}{cccc}
		\mathbf{p}_1 & \mathbf{p}_2 & ... & \mathbf{p}_{n_k} 
		\end{array}\right]\tilde{\bsym{\Pi}}\mathbf{u}^E.
	\end{split}
\end{equation}
To be clear, the strain operator acting on the row vector of polynomial base functions is
\begin{equation}\label{ep8}
	\begin{split}
		\bsym{\epsilon}\left[ \begin{array}{cccc}
		\mathbf{p}_1 & \mathbf{p}_2 & ... & \mathbf{p}_{n_k} 
		\end{array}\right]=\left[ \begin{array}{cccc}
		\partial_{x_1}p_{1,1} & \partial_{x_1}p_{2,1} & ... & \partial_{x_1}p_{n_k,1} \\
		\partial_{x_2}p_{1,2} & \partial_{x_2}p_{2,2} & ... & \partial_{x_2}p_{n_k,2} \\
		\partial_{x_2}p_{1,1}+\partial_{x_1}p_{1,2} & \partial_{x_2}p_{2,1}+\partial_{x_1}p_{2,2} & ... & \partial_{x_2}p_{n_k,1}+\partial_{x_1}p_{n_k,2} 
		\end{array}\right].
	\end{split}
\end{equation}\\

\noindent
\textbf{Remarks}
\begin{enumerate}[(i)]
	\item The strains resulting from the above work are organized using Voigt notation.  The resulting strain vector contains the two-dimensional engineering strains $\epsilon_x$, $\epsilon_y$, $\gamma_{xy}=2\epsilon_{xy}$.
	\item The strains are found for an individual polygonal (VEM) element.  Hence, they are constant within the element.
	\item In \eqref{ep8}, $p_{i,j}$ is the $j$th component of polynomial vector $\mathbf{p}_i$.  In 2D the vectors $\mathbf{p}_i$ are as indicated in \eqref{pf1} and \eqref{pf2}.
\end{enumerate}

\section{Element Stresses}
The element-wise stresses are calculated using the previously found strain vector, $\bsym{\epsilon}(\mathbf{v}^h)\approx\bsym{\epsilon}(\Pi(\mathbf{v}^h))$.  The stress vector is
\begin{equation}\label{s1}
	\bsym{\sigma}(\mathbf{v}^h)=\mathbf{C}\bsym{\epsilon}(\Pi(\mathbf{v}^h)),
\end{equation}
where for plane stress
\begin{equation}\label{s2}
	\mathbf{C}=\frac{E}{1-\nu^2}\left[ \begin{array}{ccc}
		1 & \nu & 0\\
		\nu & 1 & 0\\
		0 & 0 & \frac{1}{2}(1-\nu)
		\end{array}\right],
\end{equation}
and for plain strain
\begin{equation}\label{s3}
	\mathbf{C}=\frac{E}{(1+\nu)(1-2\nu)}\left[ \begin{array}{ccc}
		1-\nu & \nu & 0\\
		\nu & 1-\nu & 0\\
		0 & 0 & \frac{1-2\nu}{2}
		\end{array}\right].
\end{equation}\\

\noindent
\textbf{Remarks}
\begin{enumerate}[(i)]
	\item The stresses resulting from the above work are organized using Voigt notation.  The resulting stress vector contains the two-dimensional engineering stresses $\sigma_x$, $\sigma_y$, $\sigma_{xy}$.
	\item The stresses are found for an individual polygonal (VEM) element.  Hence, they are constant within the element.
	\item In \eqref{s2} and \eqref{s3}, $E$ is the modulus of elasticity and $\nu$ is Poisson's ratio.
\end{enumerate}

\section{Element Internal Forces}
With an eye toward applications with nonlinear analysis, element internal forces are calculated by multiplying the element stiffness matrix times the vector of element displacements.  For example, the internal force vector for a single element $i$ is
\begin{equation}\label{eif1}
	\mathbf{q}^i_{int}=\mathbf{k}_E\mathbf{u}^E.
\end{equation}
Then similar to FEM the individual internal force vectors for all elements are assembled into the global internal force vector using the assembly operator \cite{Hughes}.  That is,
\begin{equation}\label{eif2}
	\mathbf{F}_{int}=\overset{n_{elem}}{\underset{i=1}{\mathbf{\mathsf{A}}}}\mathbf{q}^{i}_{int}.
\end{equation}

\section{Some Relevant Concepts and Terminology}

Many concepts are used to construct VEM.  It is useful to organize and explain these concepts.  Without such an overview it is easy to get lost in the terminology and and lose sight of the ultimate objective.  The objective here is to explain concepts needed for the numerical solution of elasticity problems using VEM.  Unless evident otherwise, the definitions of variables below are for 2D.

\begin{itemize}
	\item $\Omega$, the symbol which represents the continuous domain of the 2D elasticity problem to be solved by VEM (see Figure \ref{fig1})
	\item $\Omega \subset \mathbb{R}^2$, the domain is contained in the real 2D coordinate space
	\item $\partial\Omega$, the boundary of the domain, this can be decomposed into prescribed displacement boundaries (Dirichlet or essential), $\partial\Omega_u$, and prescribed traction boundaries (Neumann or natural), $\partial\Omega_t$.  It is true that $\partial\Omega$=$\partial\Omega_u \cup \partial\Omega_t$
	\item $\mathbf{n}$, used to denote an outward normal vector to the boundary
	\item $\mathbf{n}_{e_j}$, used to denote an outward normal vector to the boundary edge $e_j$
	\item $\bar{\mathbf{u}}=\mathbf{0}$ on $\partial\Omega_u$, prescribed homogeneous displacement boundary condition
	\item $\bar{\mathbf{t}}$, prescribed traction boundary condition
	\item $\mathbf{u}$, the displacement solution to the elasticity problem.  In 2D this is just a column vector with two components ($\mathbf{u}=[u_1 \ u_2]^T$) that are functions of the coordinates $(x_1, x_2)$.
	\item $\boldsymbol{\mathcal{V}}$, defined here as a vector-valued function space, in our case in 2D, with components $v_1, v_2$.  The components belong to a first-order Sobolev space $\mathcal{H}^1(\Omega)$ with zero values on displacement boundaries.  The space contains functions that are square-integrable up to and including first derivatives.  Functions that have these characteristics are needed later.  These careful definitions help us know exactly what type of functions we want, and help us avoid problematic functions (that might give infinite square integrable results.  Such functions would imply infinite strain energy, which is not allowed.)  This function space is defined compactly as $\boldsymbol{\mathcal{V}}\equiv[\mathcal{H}_0^1(\Omega)]^2$.
	\item $a(\mathbf{u},\mathbf{v})=\int _{\Omega} \boldsymbol{\sigma}(\mathbf{u}):\boldsymbol{\epsilon}(\mathbf{v})d\Omega$, a \emph{bilinear form} related to internal strain energy used in problems of linear elasticity.
	\item $L(\mathbf{v})=\int_{\Omega} \mathbf{v}\cdot\mathbf{f}d\Omega+\int_{\partial\Omega_t} \mathbf{v}\cdot\bar{\mathbf{t}}d\partial\Omega$, a \emph{linear form} related to the external energy caused by external loads applied to the domain.  This could also include external energy caused by external prescribed displacements.  However, in this work external prescribed displacements are assumed zero for simplicity.
	\item $\boldsymbol{\sigma}=\boldsymbol{C}\boldsymbol{\epsilon}$, Hooke's law for linear elasticity relating stresses to strains.  In indicial notation this is written as $\sigma_{ij}=C_{ijkl}\epsilon_{kl}$.  In Voigt notation, in 2D, it is a 3x1 column vector.  In tensor notation, in 3D, it has 9 components and it is often expressed as a 3x3 symmetric matrix.
	\item $\boldsymbol{\epsilon}(\mathbf{u})=\frac{1}{2}(\boldsymbol{\nabla} \mathbf{u}+\boldsymbol{\nabla} \mathbf{u}^T)$, is the linearized (small) strain tensor, in indicial notation this is written as $\epsilon=\frac{1}{2}(u_{i,j}+u_{j,i})$.  In Voigt notation, in 2D, it is a 3x1 column vector.  In tensor notation, in 3D, it has 9 components and it is often expressed as a 3x3 symmetric matrix.  See also \eqref{pf7}, \eqref{ep7}, \eqref{ep8}.
	\item $\boldsymbol{\mathcal{V}}^h$, this is the discrete vector-valued function space of VEM trial solutions $\mathbf{u}^h$ and weight functions $\mathbf{v}^h$.  Exact continuous analytical solutions can sometimes (but rarely) be found for an elasticity problem.  Such solutions reside in the space of functions $\boldsymbol{\mathcal{V}}$.  However, for many problems only discrete (numerical) solutions are possible by FEM or VEM.  Hence, the domain is discretized into sub domains (elements) in which discrete functions $\mathbf{u}^h$are used to approximate the elementwise solution.  These functions are piecewise connected, at element boundaries, across the problem domain.  The discrete space of functions is a subset of the space that includes continuous analytical functions ($\boldsymbol{\mathcal{V}}^h$ $\ \subset \ $ $\boldsymbol{\mathcal{V}}$).  This space of functions $\boldsymbol{\mathcal{V}}(E)^h$ on elements $E$ contains polynomial functions as well as non-polynomial functions.  For a more formal definition, see Appendix A.1 of Mengolini et al. \cite{mengolini}.
	\item $E$, an individual polygon domain
	\item $\Omega^h$, the discretized domain, covered by a collection of elements
	\item VEM functions, the functions used in the virtual element method are found in the space of functions $\boldsymbol{\mathcal{V}}^h$.  VEM functions include a combination of polynomial and non-polynomial type functions.
	\item $\bsym{\varphi}$, VEM shape function matrix, a $2\times2n_d$ matrix
	\item $\bsym{\varphi}_i$, VEM shape function vector associated with element degree of freedom $i$, a $2\times 1$ column vector
	\item $\boldsymbol{\mathcal{P}}_k(E)$, the space of polynomial functions of order less than or equal to $k$
	\item $\Pi_{E,k}$, the local projection operator.  This operator projects VEM functions onto the space of polynomials of order $k$ or less.  In math terms this is expressed as $\Pi_{E,k}$ \ : $\boldsymbol{\mathcal{V}}(E)^h$ $\ \rightarrow \ $ $\boldsymbol{\mathcal{P}}_k(E)$
	\item $\Pi$, the projector operator, to be understood as a simplified version of $\Pi_{E,k}$, unless directed otherwise
	\item $n_d$, number of degrees of freedom along one spatial dimension of an element
	\item $n_v$, number of polygon vertices for element $E$.  Importantly, for $k=1$ the number of degrees of freedom along one spatial dimension, $n_d$ equals the number of vertices, $n_v$.
	\item $n_k$, the number of terms in in the polynomial base of order $k$.  That is, $n_k=(k+1)(k+2)$
	\item $\bsym{\partial}$, the differential operator defined in \eqref{pf3}, a $3\times2$ operator matrix
	\item $dof_i(\mathbf{v}^h)$, degree of freedom $i$ of $\mathbf{v}^h$ for element $E$.	
	\item $\mathbf{u}^h$, discrete trial solution, a $2\times 1$ column vector
	\item $\mathbf{v}^h$, discrete weight function, a $2\times 1$ column vector
	\item $\Delta\mathbf{v}^h|_E=\Delta^2\mathbf{v}^h|_E=$ Laplacian of $\mathbf{v}^h|_E$ over element $E$
	\item $k$, degree of polynomials used to approximate displacement within each element, degree of the polynomial base
	\item $e_j$, element edge $j$
	\item $|e_j|$, length of element edge $j$
	\item $|E|$, area of element $E$
	\item $\mathbf{B}$, strain displacement operator acting on VEM shape functions, a $3\times2n_d$ matrix
	\item $\bar{\mathbf{B}}$, the final modified $n_k\times2n_d$ ``B" matrix used to calculate the projector, $\tilde{\bsym{\Pi}}=\mathbf{G}^{-1}\bar{\mathbf{B}}$
	\item $\tilde{\mathbf{B}}$, the $n_k\times2n_d$ ``B" matrix that results in an undetermined system for the projector, $\tilde{\mathbf{B}}=\tilde{\mathbf{G}}\tilde{\bsym{\Pi}}^*$, it is the ``B" matrix that needs its first three rows modified with $\breve{\mathbf{B}}$ to get $\bar{\mathbf{B}}$.
	\item $\breve{\mathbf{B}}$, this is the $3\times2n_d$ matrix that is inserted into the first three rows of $\tilde{\mathbf{B}}$ to get the final modified matrix $\bar{\mathbf{B}}$.
	\item $\mathbf{G}$, the final $n_k\times n_k$ ``G" matrix used to calculate the projector, $\tilde{\bsym{\Pi}}=\mathbf{G}^{-1}\bar{\mathbf{B}}$
	\item $\tilde{\mathbf{G}}$, the $n_k\times n_k$ ``G" matrix that is part of the undetermined system $\tilde{\mathbf{B}}=\tilde{\mathbf{G}}\tilde{\bsym{\Pi}}^*$
	\item $\breve{\mathbf{G}}$, the $3\times n_k$ matrix that is inserted into the first three rows of $\tilde{\mathbf{G}}$ to obtain $\mathbf{G}$
	\item $\bsym{\Pi}$, the $2n_d\times2n_d$ projector matrix that is used to construct the stability stiffness, $\mathbf{k}^s_E$.  It is the energy projector operator with respect to the $\bsym{\varphi}$ basis set \cite{suku:elastodyn}. 
	\item $\tilde{\bsym{\Pi}}$, the $n_k\times2n_d$ projector matrix that is used to construct the consistency stiffness, $\mathbf{k}^c_E$.  It is the energy projector operator with respect to the polynomial basis set \cite{suku:elastodyn}.
	\item $\tilde{\bsym{\Pi}}^*$, the $n_k\times2n_d$ projector matrix that is part of the undetermined system, $\tilde{\mathbf{B}}=\tilde{\mathbf{G}}\tilde{\bsym{\Pi}}^*$
	\item $\breve{\bsym{\Pi}}$, the $3\times2n_d$ matrix that relates $\breve{\mathbf{G}}$ and $\breve{\mathbf{B}}$
	\item $\mathbf{D}$, this $2n_d\times n_k$ matrix is ``used to express the projection of a VEM function as a linear combination of the VEM functions themselves" \cite{mengolini}.
	\item $\mathbf{P}_1$, the $2\times n_k$ polynomial basis matrix of order $k=1$
	\item $\mathbf{p}_i$, an individual scaled vector monomial in the polynomial basis set, a $2\times1$ column vector
	\item $\mathbf{k}^c_E$, the $2n_d\times2n_d$ stiffness matrix for element $E$ that provides consistency
	\item $\mathbf{k}^s_E$, the $2n_d\times2n_d$ stiffness matrix for element $E$ that provides stability	
\end{itemize}

	\begin{figure}
  \centering
  \epsfig{file=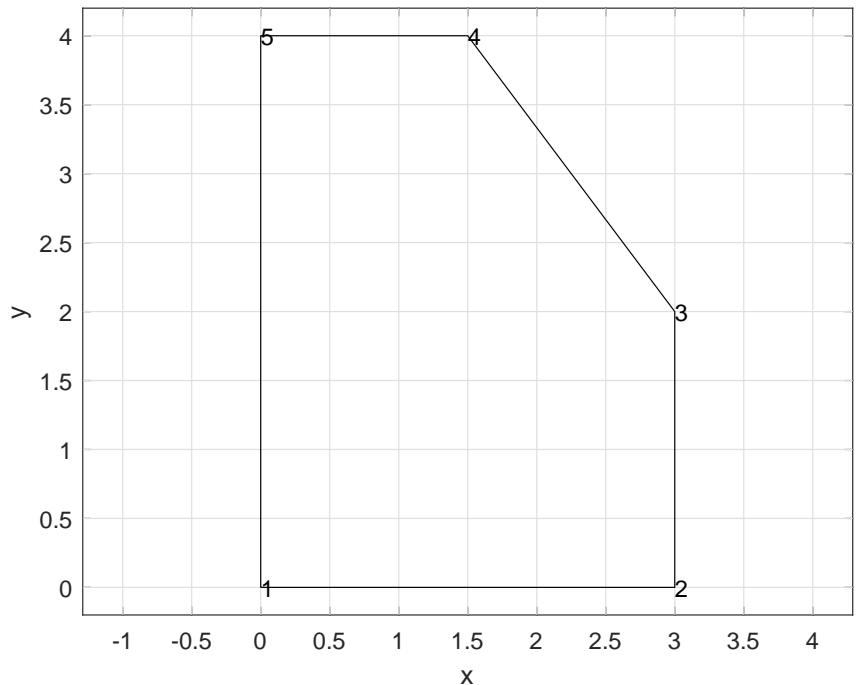,width=.8\textwidth}
  \caption{Single five sided element with vertex labels shown.}\label{fig3}
  \end{figure}  

\section{Example -- A single 5 sided element}
Various terms in the VEM formulation are calculated for a single 5 sided element.  The numerical results are provided so that readers implementing VEM can verify that calculations are correct.  Element geometry is provided in figure \ref{fig3}.\\

Input
\begin{itemize}
	\item Polynomial degree on polygon edges: $k=1$
	\item Modulus of Elasticity: $E=1000$
	\item Poisson's ratio: $\nu=0.3$
	\item Plane stress problem
	\item Element Node Numbers: [1,2,3,4,5]
	\item Domain Thickness: t=1
	\item Specified zero displacements: Node 1 $(u_x=0,u_y=0)$, Node 5 $(u_x=0)$	
	\item Specified point loads: Node 2 $(F_x=40)$, Node 3 $(F_x=80)$, Node 4 $(F_x=40)$
\end{itemize}
Output
\begin{itemize}
	\item Centroid location: $\bar{x}=1.3571, \ \bar{y}=1.8095$
	\item Number of vertices: $n_v=5$
	\item Number of dofs: $2n_d=10$
	\item Polygon Diameter: $h_E=5$
	\item Polygon Area: $|E|=10.5$
	\item $\bar{\mathbf{B}}$ matrix:
	\begingroup
    \fontsize{8pt}{12pt}\selectfont
	\begin{verbatim}
    0.2000         0    0.2000         0    0.2000         0    0.2000         0    0.2000         0
         0    0.2000         0    0.2000         0    0.2000         0    0.2000         0    0.2000
    0.0724   -0.0543    0.0724    0.0657   -0.0076    0.0657   -0.0876    0.0057   -0.0876   -0.0543
 -230.7692 -307.6923 -230.7692  153.8462  115.3846  307.6923  230.7692  153.8462  115.3846 -307.6923
 -439.5604  -98.9011  219.7802  -98.9011  439.5604   49.4505  219.7802   98.9011 -439.5604   49.4505
 -131.8681 -329.6703   65.9341 -329.6703  131.8681  164.8352   65.9341  329.6703 -131.8681  164.8352		
	\end{verbatim}
	\endgroup
	\item $\mathbf{D}$ matrix:
	\begingroup
    \fontsize{8pt}{12pt}\selectfont
	\begin{verbatim}
    1.0000         0    0.3619   -0.3619   -0.2714         0
         0    1.0000   -0.2714   -0.2714         0   -0.3619
    1.0000         0    0.3619   -0.3619    0.3286         0
         0    1.0000    0.3286    0.3286         0   -0.3619
    1.0000         0   -0.0381    0.0381    0.3286         0
         0    1.0000    0.3286    0.3286         0    0.0381
    1.0000         0   -0.4381    0.4381    0.0286         0
         0    1.0000    0.0286    0.0286         0    0.4381
    1.0000         0   -0.4381    0.4381   -0.2714         0
         0    1.0000   -0.2714   -0.2714         0    0.4381	
	\end{verbatim}
	\endgroup	
	\item $\mathbf{G}$ matrix:
	\begingroup
    \fontsize{8pt}{12pt}\selectfont
	\begin{verbatim}
    1.0000         0   -0.0381    0.0381    0.0286         0
         0    1.0000    0.0286    0.0286         0    0.0381
   -0.0381    0.0286    0.2023   -0.0566    0.0229   -0.0229
    0.0000         0   -0.0000  646.1538         0         0
         0         0   -0.0000   -0.0000  461.5385  138.4615
    0.0000         0    0.0000    0.0000  138.4615  461.5385
	\end{verbatim}
	\endgroup	
	\item $\tilde{\bsym{\Pi}}$ matrix:
	\begingroup
    \fontsize{8pt}{12pt}\selectfont
	\begin{verbatim}
    0.2566   -0.0016    0.2093    0.0016    0.1635   -0.0000    0.1592   -0.0033    0.2114    0.0033
   -0.0016    0.2556    0.0033    0.2124   -0.0033    0.1592    0.0000    0.1616    0.0016    0.2112
    0.4143   -0.5190    0.2429    0.2810   -0.0643    0.4762   -0.3571    0.1524   -0.2357   -0.3905
   -0.3571   -0.4762   -0.3571    0.2381    0.1786    0.4762    0.3571    0.2381    0.1786   -0.4762
   -0.9524         0    0.4762    0.0000    0.9524    0.0000    0.4762         0   -0.9524   -0.0000
         0   -0.7143    0.0000   -0.7143    0.0000    0.3571   -0.0000    0.7143         0    0.3571
	\end{verbatim}
	\endgroup	
	\item $\bsym{\Pi}$ matrix:
	\begingroup
    \fontsize{8pt}{12pt}\selectfont
	\begin{verbatim}
    0.7943   -0.0171    0.2971    0.0171   -0.1829   -0.0000   -0.2286   -0.0343    0.3200    0.0343
   -0.0171    0.7843    0.0343    0.3300   -0.0343   -0.2286    0.0000   -0.2029    0.0171    0.3171
    0.2229   -0.0171    0.5829    0.0171    0.3886    0.0000    0.0571   -0.0343   -0.2514    0.0343
    0.0171    0.1871   -0.0343    0.6414    0.0343    0.3429   -0.0000    0.0314   -0.0171   -0.2029
   -0.0857   -0.0000    0.3429    0.0000    0.4857    0.0000    0.3429   -0.0000   -0.0857   -0.0000
    0.0171   -0.0986   -0.0343    0.3557    0.0343    0.4857   -0.0000    0.3171   -0.0171   -0.0600
   -0.1086    0.0171   -0.0400   -0.0171    0.2971    0.0000    0.4857    0.0343    0.3657   -0.0343
    0.0000   -0.0857    0.0000   -0.0857    0.0000    0.3429   -0.0000    0.4857   -0.0000    0.3429
    0.1771    0.0171   -0.1829   -0.0171    0.0114   -0.0000    0.3429    0.0343    0.6514   -0.0343
   -0.0171    0.2129    0.0343   -0.2414   -0.0343    0.0571   -0.0000    0.3686    0.0171    0.6029
	\end{verbatim}
	\endgroup	
	\item $\mathbf{k}_E$, element stiffness matrix:
	\begingroup
    \fontsize{8pt}{12pt}\selectfont
	\begin{verbatim}
  523.2489  204.4601 -159.9480   38.8680 -438.1401 -156.9859 -269.0252 -148.3797  343.8645   62.0375
  204.4601  404.4220   62.0375  128.4422 -148.3797 -241.5527 -156.9859 -286.5997   38.8680   -4.7119
 -159.9480   62.0375  251.9156 -101.2839  104.5264  -86.3422   19.7167   -9.3631 -216.2107  134.9518
   38.8680  128.4422 -101.2839  338.6842  -67.4759 -110.0770    7.8493 -200.8041  122.0425 -156.2453
 -438.1401 -148.3797  104.5264  -67.4759  522.9966  102.0408  210.1555  123.1778 -399.5384   -9.3631
 -156.9859 -241.5527  -86.3422 -110.0770  102.0408  291.1714  133.4380  150.6317    7.8493  -90.1734
 -269.0252 -156.9859   19.7167    7.8493  210.1555  133.4380  272.8564  102.0408 -233.7034  -86.3422
 -148.3797 -286.5997   -9.3631 -200.8041  123.1778  150.6317  102.0408  356.7551  -67.4759  -19.9830
  343.8645   38.8680 -216.2107  122.0425 -399.5384    7.8493 -233.7034  -67.4759  505.5879 -101.2839
   62.0375   -4.7119  134.9518 -156.2453   -9.3631  -90.1734  -86.3422  -19.9830 -101.2839  271.1137
	\end{verbatim}
	\endgroup
	\item $u_x, \ u_y$, nodal displacements:
	\begingroup
    \fontsize{8pt}{12pt}\selectfont
	\begin{verbatim}
0.00  0.000
0.12  0.000
0.12 -0.024
0.06 -0.048
0.00 -0.048
	\end{verbatim}
	\endgroup			
	
	\item $\bsym{\epsilon}$, strains ($\epsilon_x$, $\epsilon_y$, $\gamma_{xy}=2\epsilon_{xy}$):
	\begingroup
    \fontsize{8pt}{12pt}\selectfont
	\begin{verbatim}
    0.0400
   -0.0120
   -0.0000
	\end{verbatim}
	\endgroup
	\item $\bsym{\sigma}$, stresses ($\sigma_x$, $\sigma_y$, $\sigma_{xy})$:
	\begingroup
    \fontsize{8pt}{12pt}\selectfont
	\begin{verbatim}
   40.0000
   -0.0000
   -0.0000
	\end{verbatim}
	\endgroup		
\end{itemize}

\section{Example -- Cantilever}
\begin{figure}
\centering
\subfigure[]{\epsfig{file =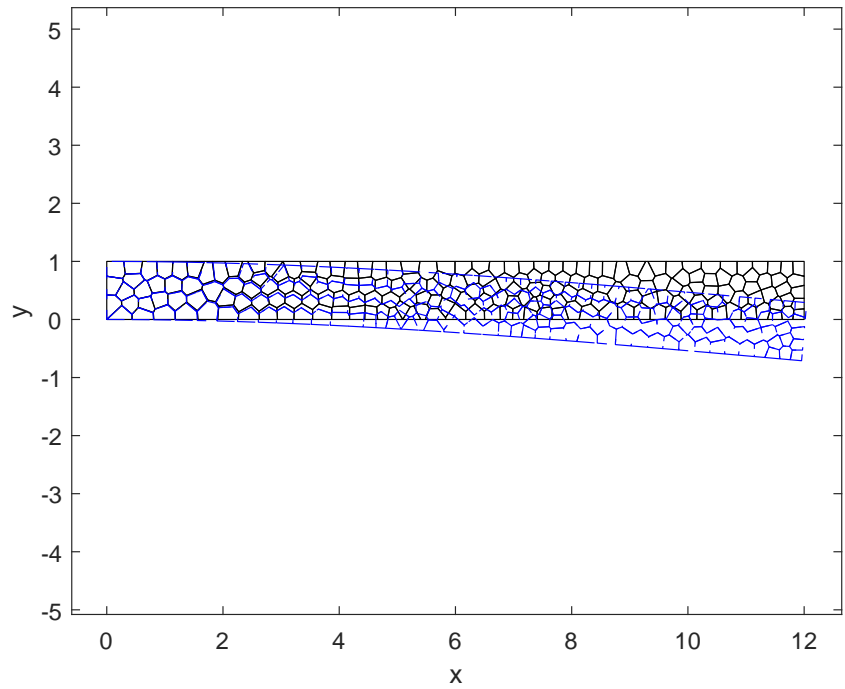,width=0.8\textwidth}}
\subfigure[]{\epsfig{file = 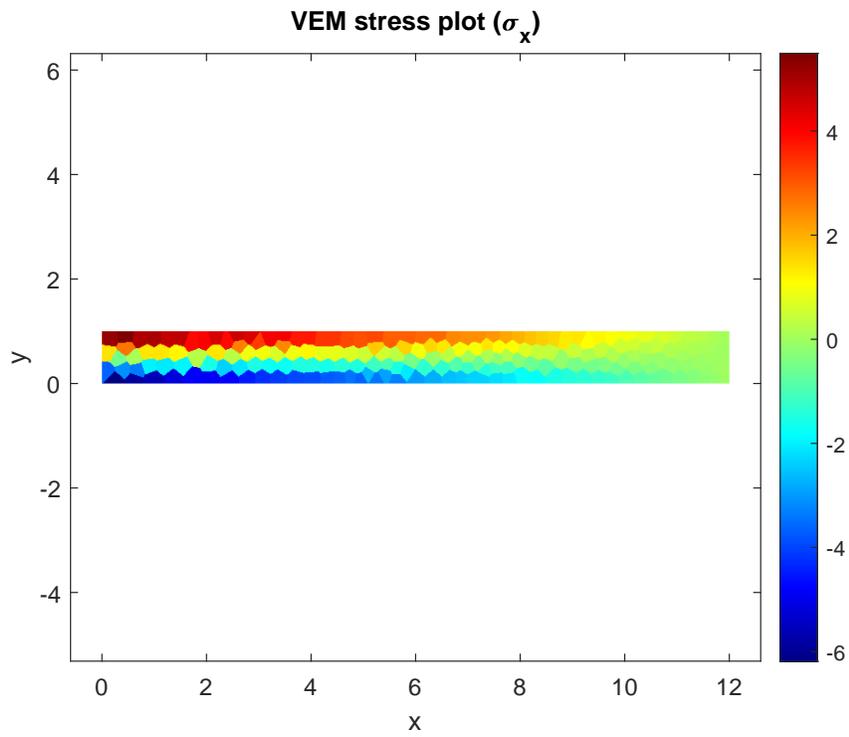,width=0.8\textwidth}}
\caption{End Loaded Cantilever: (a)~Original and deformed shape, (b)~Bending stresses,~$\sigma_x$.}
\label{fig4}
\end{figure}

A 12 inch long cantilever is loaded with a point load at its free end.  The cantilever is 1 inch deep and 1 inch thick into the page.  The load at the end is 0.1 kips.  The modulus of elasticity is 1000 ksi and Poisson's ratio is 0.3.  The deflected shape is shown in Figure \ref{fig4}a.  The bending stresses, $\sigma_x$, are shown in Figure \ref{fig3}b.  It is evident that stresses are constant over each polygon element according to the VEM formulation.  The cantilever has all nodes pinned in the x and y direction at the support for this example.  The maximum bending stress is 6.19 ksi compared to the theoretical prediction of 7.2 ksi.  A finer discretization of the domain would provided better results.  In this example, polymesher \cite{talischi} was used to randomly discretize the domain with 200 polygons.  The tip displacement for this example is 0.71 inches and the predicted value is 0.691 inches, according to the simple beam theory formula, $\Delta=\frac{PL^3}{3EI}$.

\section{Example -- Plate with hole}
\begin{figure}
\centering
\subfigure[]{\epsfig{file =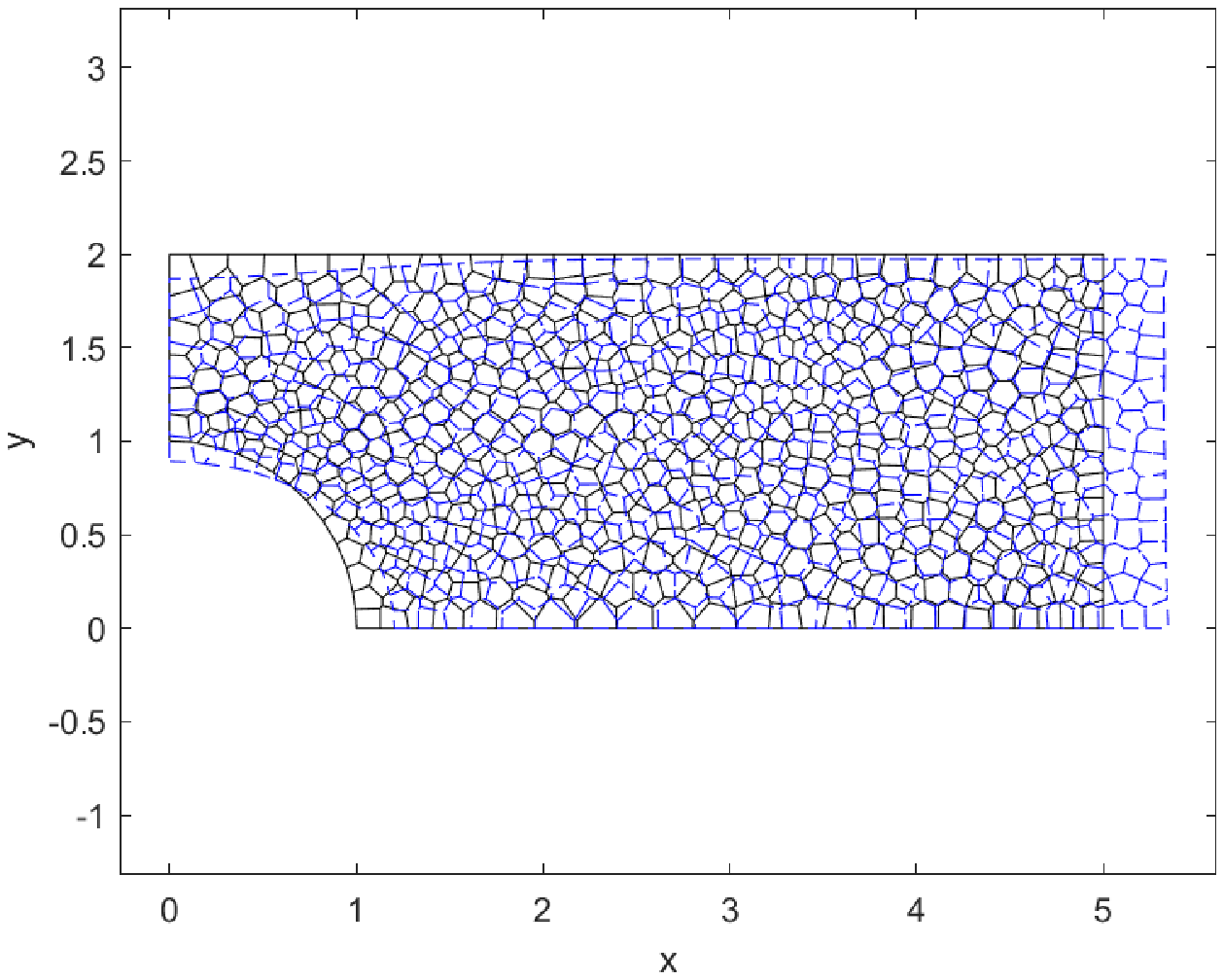,width=0.5\textwidth}}
\subfigure[]{\epsfig{file = 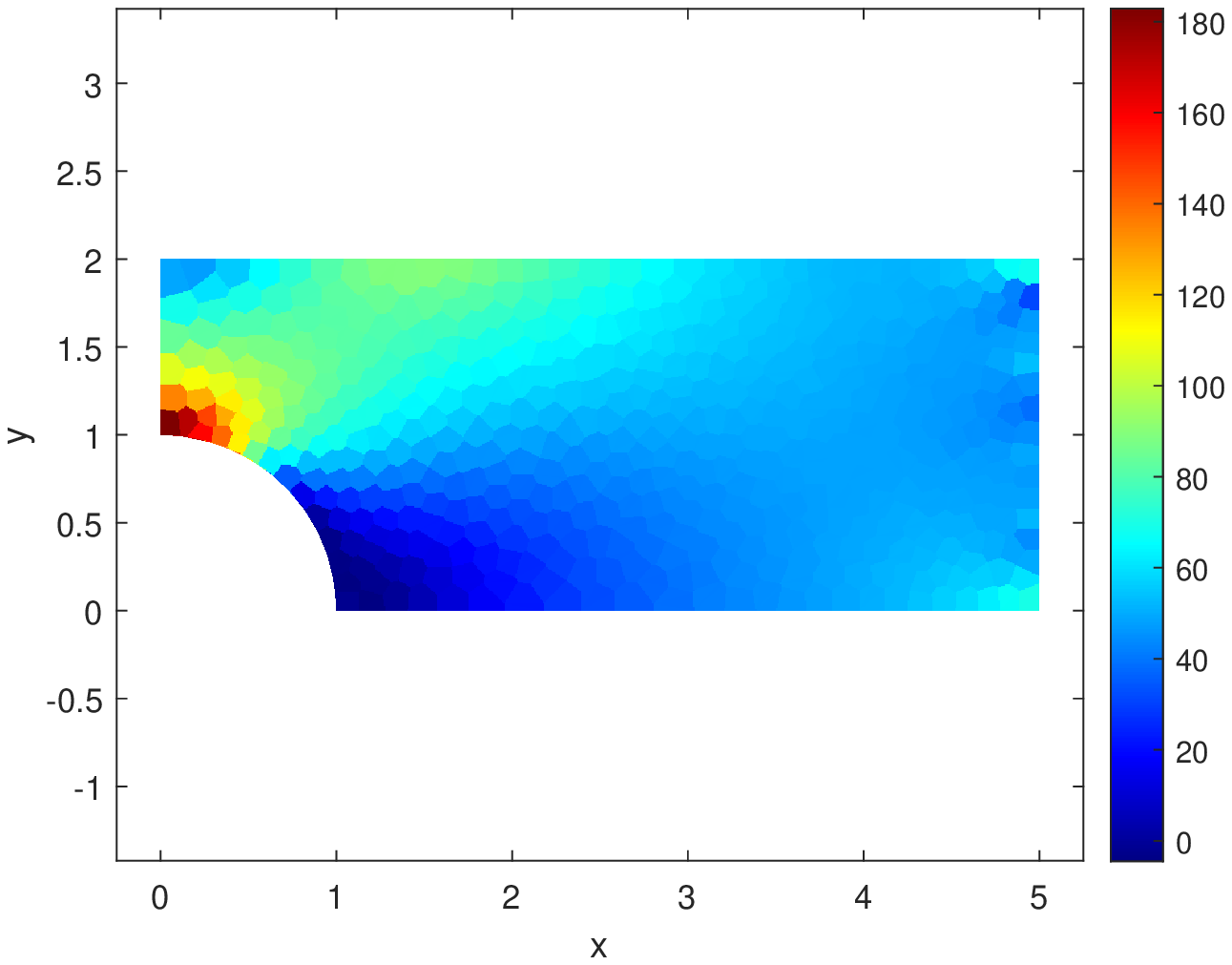,width=0.5\textwidth}}
\subfigure[]{\epsfig{file = 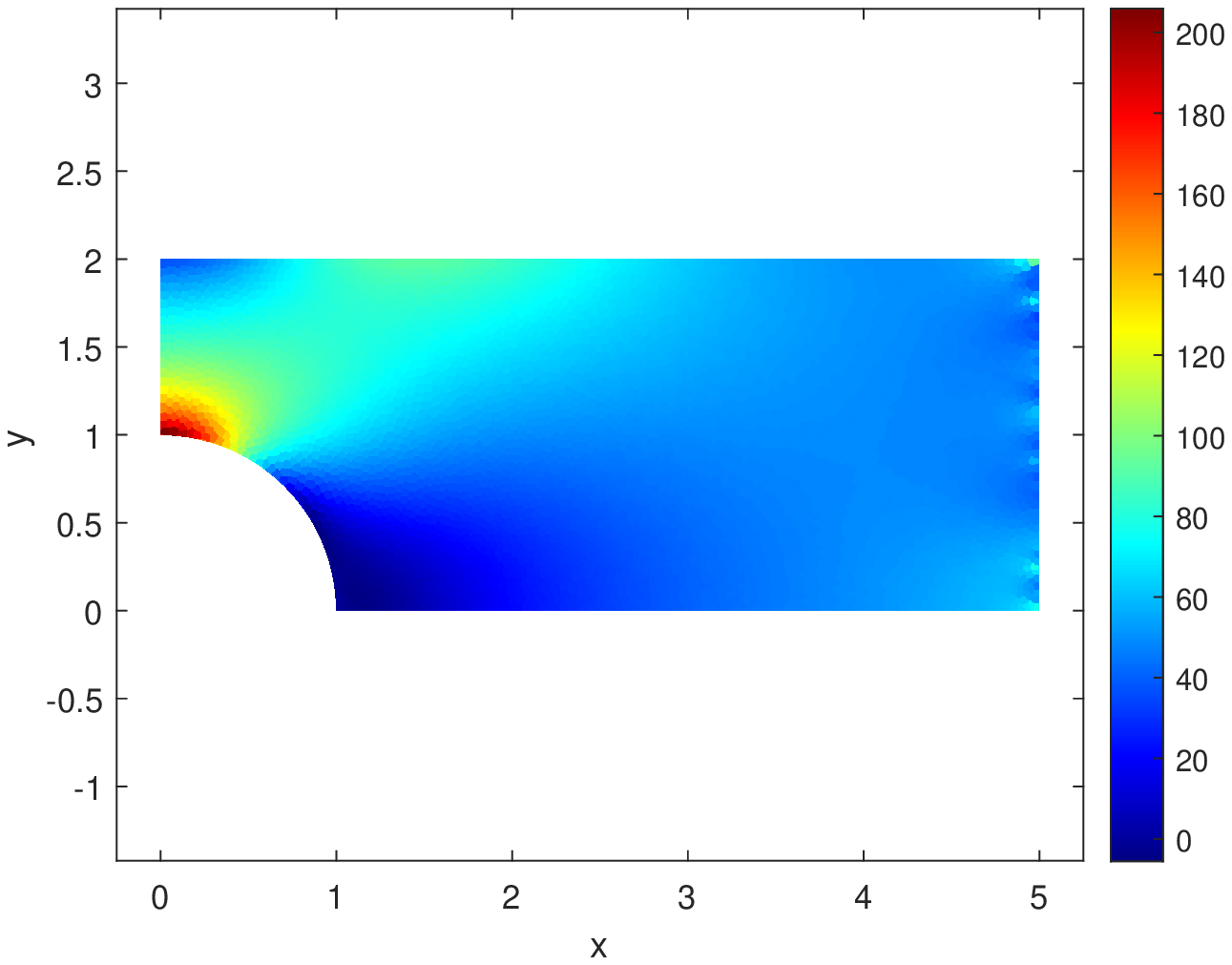,width=0.5\textwidth}}
\caption{Quadrant of plate with hole in tension: (a)~Original and deformed shape, (b)~$\sigma_x$ stresses with 500 polygon discretization, (c)~$\sigma_x$ stresses with 5000 polygon discretization.}
\label{fig5}
\end{figure}

A plate with a hole is loaded in tension.  Due to symmetry only one quadrant of the plate is analyzed. The modulus of elasticity is 1000 ksi and Poisson's ratio is 0.3.  The deflected shape is shown in Figure \ref{fig5}a.  The $\sigma_x$ stresses are shown in Figures \ref{fig5}b,c with 500 and 5000 polygons respectively.  It is evident that stresses are constant over each polygon element according to the VEM formulation.  The plate has zero displacement supports in the x-direction at the left vertical edge, and y-direction at the bottom horizontal edge for this example.

\section{Conclusion}
An explanation of the 2D virtual element method (VEM) is provided.  Detailed derivations and numerical examples are given.  It is shown that VEM is a viable alternative to standard FEM formulations.

\section{Acknowledgments}
The author would like to thank Professor N. Sukumar for helpful discussions regarding the Virtual Element Method.  Yet, any errors or conceptual shortcomings are entirely the responsibility of the author.

\bibliographystyle{hsiam}
\bibliography{IntroToVEM}
\end{document}